\documentclass[12pt,leqno]{amsart} 
\newtheorem{thm}{Theorem}[section]
\newtheorem{defi}{Definition}[section]
\newtheorem{lem}{Lemma}[section]
\newtheorem{cor}{Corollary}[section]
\newtheorem{pr}{Proposition}[section]
\newtheorem{rem}{Remark}[section]
\newcommand{\be}{\begin{equation}}
\newcommand{\ee}{\end{equation}}
\newcommand{\bea}{\begin{eqnarray}}
\newcommand{\eea}{\end{eqnarray}}
\newcommand{\beb}{\begin{eqnarray*}}
\newcommand{\eeb}{\end{eqnarray*}}
\usepackage{amssymb,amsfonts,amsthm,setspace,indentfirst,mathrsfs}
\usepackage{minibox}
\numberwithin{equation}{section}
\begin{document}
%%%%%%%%%%%%%%%%%%%%%%%%%%%%%%%%%%%%%%%%%%%%%%%%%%%%%%%%%%%%%%%%%%%%%%%%%%%%%%%%%%%%%%%%%%%%%%%%%%%%%%%%%%%%%
%
\title[On some curvature restricted geometric structures]
{\bf{On some curvature restricted geometric structures for projective curvature tensor}}
\author[Absos Ali Shaikh and Haradhan Kundu]{Absos Ali Shaikh and Haradhan Kundu}
\date{\today}
\address{\noindent\newline Department of Mathematics,\newline The University of 
Burdwan, Golapbag,\newline Burdwan-713104,\newline West Bengal, India}
\email{aask2003@yahoo.co.in, aashaikh@math.buruniv.ac.in}
\email{kundu.haradhan@gmail.com}
\dedicatory{}
%%%%%%%%%%%%%%%%%%%%%%%%%%%%%%%%%%%%%%%%%%%%% Abstract %%%%%%%%%%%%%%%%%%%%%%%%%%%%%%%%%%%%%%%%%%%%%%%%%%%%%
\begin{abstract}
The projective curvature tensor $P$ is invariant under a geodesic preserving transformation on a semi-Riemannian manifold. It is well known that $P$ is not a generalized curvature tensor and hence it possesses different geometric properties than other generalized curvature tensors. The main object of the present paper is to study some semisymmetric type and pseudosymmetric type curvature restricted geometric structures due to projective curvature tensor. The reduced pseudosymmetric type structures for various Walker type conditions are deduced and the existence of Venzi space is ensured. It is shown that the geometric structures formed by imposing projective operator on a (0,4)-tensor is different from that for the corresponding (1,3)-tensor. Characterization of various semisymmetric type and pseudosymmetric type curvature restricted geometric structures due to projective curvature tensor is obtained on a Riemannian and a semi-Riemannian manifold, and it is shown that some of them reduce to Einstein manifold for the Riemann case. Finally to support our theorems four suitable examples are presented.
\end{abstract}
%%%%%%%%%%%%%%%%%%%%%%%%%%%%%%%%%%%%%%%%%%%%%%%%%%%%%%%%%%%%%%%%%%%%%%%%%%%%%%%%%%%%%%%%%%%%%%%%%%%%%%%%%%%%
%
\maketitle
%

%%%%%%%%%%%%%%%%%%%%%%%%%%%%%%%%%%%%%%%%%%%%%%%%%%%%%%%%%%%%%%%%%%%%%%%%%%%%%%%%%%%%%%%%%%%%%%%%%%%%%%%%%%%%%%
%																					Introduction
%%%%%%%%%%%%%%%%%%%%%%%%%%%%%%%%%%%%%%%%%%%%%%%%%%%%%%%%%%%%%%%%%%%%%%%%%%%%%%%%%%%%%%%%%%%%%%%%%%%%%%%%%%%%%%
\section{Introduction}
%==========================
Let $M$ be an $n$-dimensional $(n\geq 3)$ connected smooth manifold equipped with the semi-Riemannian metric $g$, Levi-Civita connection $\nabla$, Riemann-Christoffel curvature tensor $R$ of type $(0,4)$, Riemann-Christoffel curvature tensor $\mathcal R$ of type $(1,3)$, Ricci tensor $S$ of type $(0,2)$, Ricci operator $\mathcal S$, scalar curvature $\kappa$, Gaussian curvature $G$ and concircular curvature tensor $W$.\\
\indent Symmetry plays an important role in the study of differential geometry of manifolds. The manifold $M$ is said to be locally symmetric (\cite{Cart26}, \cite{Cart27}, \cite{Cart46}) if its local geodesic symmetries are isometry and $M$ is said to be globally symmetric if its geodesic symmetries are extendible to the whole of $M$. In terms of curvature restriction $M$ is locally symmetric if $\nabla R = 0$ (see \cite{Cart46} and also \cite{CP70}, \cite{CP80}). We note that a geometric structure on $M$ formed by imposing a restriction on some curvature tensors of $M$ is called a curvature restricted geometric structure. During the last eight decades the notion of local symmetry has been generalized by many authors by weakening the restriction $\nabla R = 0$ and there arose various curvature restricted geometric structures.\\
\indent Generalizing the notion of local symmetry, Cartan \cite{Cart46} (see also \cite{Szab82}, \cite{Szab84}, \cite{Szab85}) introduced the notion of semisymmetric manifold. A semi-Riemannian manifold $M$ is said to be semisymmetric \cite{Cart46} if 
$$\mathscr R(X,Y)\cdot R = 0,$$
where $X, Y \in \chi(M)$, the set of all smooth vector fields on $M$ and $\mathscr R(X,Y)$ is the curvature operator corresponding to $R$. Again during the study of totally umbilical submanifolds of semisymmetric manifolds as well as during the consideration of geodesic mappings on semisymmetric manifolds, Adamow and Deszcz \cite{AD83} (see also \cite{Desz92} and references therein) introduced the notion of pseudosymmetric manifolds which generalizes the notion of semisymmetric manifolds. A semi-Riemannian manifold $M$ is said to be pseudosymmetric if  
$$
\mbox{$\mathcal R(X,Y)\cdot R$ and $\mathcal G(X,Y)\cdot R$ are linearly dependent,}
$$
where $\mathcal G (X,Y) = X\wedge Y$ is the curvature operator corresponding to the Gaussian curvature tensor $G$. Replacing $\mathcal R(X,Y)$, $\mathcal G(X,Y)$ and $R$ by other curvature tensors in the defining condition of semisymmetric manifold and pseudosymmetric manifold one can get various curvature restricted geometric structures, which are respectively known as semisymmetric type and pseudosymmetric type manifolds. Deszcz and his coauthors (see \cite{DK03}, \cite{DH03}, \cite{DGHZ15} and also references therein) studied various pseudosymmetric type curvature restricted geometric structures. Recently, the present authors \cite{SK14} classified various curvature restricted geometric structures (especially, semisymmetric and pseudosymmetric) and studied their equivalency.\\
\indent The geodesic preserving transformation between two semi-Riemannian manifolds is called projective transformation and the projective curvature tensor $P$, given by
$$P(X_1,X_2,X_3,X_4) = R(X_1,X_2,X_3,X_4) - \frac{1}{n-1} \left[S(X_2,X_3)g(X_1,X_4) - S(X_1,X_3)g(X_2,X_4)\right],$$
is an invariant under such a transformation, $X_i \in \chi(M)$. A $(0,4)$ tensor is called a generalized curvature tensor if it obeys the symmetries like $R$. We note that $P$ is not a generalized curvature tensor since
$$P(X_1,X_2,X_3,X_4) \neq P(X_3,X_4,X_1,X_2) \ \mbox{in general.}$$
\indent The main object of the present paper is to study various semisymmetric type and pseudosymmetric type curvature restricted geometric structures due to the projective curvature tensor. Since $P$ is not a generalized curvature tensor, the projective curvature operator $\mathscr P(X,Y)$ can not commute with contraction. As a consequence the structures formed by the curvature operator $\mathscr P(X,Y)$ gives various interesting results different from the other generalized curvature tensors. For example the geometric structures formed by imposing $\mathcal P(X,Y)$ to a (0,4)-tensor $H$ and to the corresponding (1,3)-tensor $\mathcal H$ are different. Various curvature restricted geometric structures for $P$ along with some additional assumptions were studied by many authors (see \cite{DS10}, \cite{SaP15}, \cite{YD10}) but they do not mention the above interesting geometric fact.\\
The main results of the paper are highlighted below:\\
(1) Established some Walker type identities and found out the necessary and sufficient conditions of various Walker type conditions formed by $P$.\\
(2) Characterized the $P$-space by Venzi and showed that such a space is of constant curvature in Riemann case and for the semi-Riemann case such a space satisfies $W\cdot W = Q(S-\frac{\kappa}{n}g, W)$.\\
(3) Showed that the geometric structures formed by applying $\mathscr P(X,Y)$ on a (0,4) tensor and the corresponding (1,3)-tensor are different and also found out the sufficient condition of their equivalency.\\
(4) Characterized the semi-Riemannian manifold satisfying the following semisymmetric and pseudosymmetric type curvature conditions:\\
(i) $P\cdot \mathcal S = 0$, (ii) $P\cdot \mathcal S = L Q(g,\mathcal S)$, (iii) $P\cdot \mathcal R = 0$, (iv) $P\cdot \mathcal R = L Q(g,\mathcal R)$, (v) $P\cdot \mathcal R = L Q(S,\mathcal R)$, (vi) $P\cdot P = 0$, (vii) $P\cdot P = L Q(g,P)$, (viii) $P\cdot P = L Q(S,P)$ (ix) $P\cdot \mathcal P = 0$, (x) $P\cdot \mathcal P = L Q(g,\mathcal P)$, (xi) $P\cdot \mathcal P = L Q(S,\mathcal P)$.\\
(5) Mentioned various curvature restricted geometric structures which are properly exist for semi-Riemannian case but in Riemann space they become Einstein.\\
(6) Showed that on a generalized Roter type Riemannian manifold, various curvature restricted geometric structures, such as $P\cdot \mathcal S = 0,$ $P\cdot \mathcal R = 0,$ $P\cdot \mathcal R = L Q(g,\mathcal R)$, $P\cdot P = L Q(g,P)$  etc. are equivalent to the manifold of constant curvature, which generalizes the main results of \cite{DRV89}.\\
\indent The paper is organized as follows. Section 2 deals with preliminaries. Section 3 is concerned with some curvature related properties. In Section 4 we present our main results. Finally, in Section 5 we present some examples to support our results.
%
%%%%%%%%%%%%%%%%%%%%%%%%%%%%%%%%%%%%%%%%%%%%%%%%%%%%%%%%%%%%%%%%%%%%%%%%%%%%%%%%%%%%%%%%%%%%%%%%%%%%%%%%%%%%%%%%%%%%
%                                               Preliminaries
%%%%%%%%%%%%%%%%%%%%%%%%%%%%%%%%%%%%%%%%%%%%%%%%%%%%%%%%%%%%%%%%%%%%%%%%%%%%%%%%%%%%%%%%%%%%%%%%%%%%%%%%%%%%%%%%%%%%
\section{Preliminaries}
%=====================
Let us consider the following notations related to $(M,g)$:\\
$C^{\infty}(M) =$ the algebra of all smooth functions on $M$, \\
$\chi(M) =$ the Lie algebra of all smooth vector fields on $M$, \\
$\chi^*(M) =$ the Lie algebra of all smooth 1-forms on $M$,\\
$\Xi(M) =$ the space of all endomorphisms on $\chi(M)$ and \\
$\mathcal T^r_k(M) =$ the space of all smooth tensor fields of type $(r,k)$ on $M$.\\
%-------------------------------------
\indent For $A, E \in \mathcal T^0_2(M)$, their Kulkarni-Nomizu product (\cite{DGHS11}, \cite{Glog02}) $A\wedge E \in \mathcal T^0_4(M)$ is given by
\beb
(A \wedge E)(X_1,X_2,X_3,X_4)&=&A(X_1,X_4)E(X_2,X_3) + A(X_2,X_3)E(X_1,X_4)\\\nonumber
&-&A(X_1,X_3)E(X_2,X_4) - A(X_2,X_4)E(X_1,X_3),
\eeb
where $X_1,X_2,X_3,X_4 \in \chi(M)$. Throughout the paper we consider $X,Y,X_1,X_2,\cdots \in \chi(M)$.\\
%----------------------------------------
\indent Again for a symmetric $(0,2)$ tensor $A$ and $X, Y \in \chi(M)$, we get $(X\wedge_A Y), \mathcal A, \mathcal A^2 \in \Xi(M)$ and $A^2\in \mathcal T^0_2(M)$, (\cite{SKgrt}, \cite{SKgrtw}) defined as follows:
$$(X\wedge_A Y)X_1 = A(Y,X_1)X-A(X,X_1)Y, \ \ g(\mathcal A X,Y) = A(X,Y),$$
$$\mathcal A^2 = \mathcal A \circ \mathcal A \ \ \mbox{and} \ \ A^2(X,Y) = g(\mathcal A^2 X, Y).$$
%-----------------------------------------------
Now  for $A \in \mathcal T^0_2(M)$ and $H \in \mathcal T^0_k(M)$, $k\ge 2$, one can define $A \wedge H$ and $X\wedge_H Y$ (\cite{ADEHM14}, \cite{SKgrt}) as follows:
\beb
(A \wedge H)(X_1,X_2,Y_1,Y_2,\cdots,Y_k)&=&A(X_1,Y_2)H(X_2,Y_1,\cdots,Y_k) + A(X_2,Y_1)H(X_1,Y_2,\cdots,Y_k)\\
&-&A(X_1,Y_1)H(X_2,Y_2,\cdots,Y_k) - A(X_2,Y_2)H(X_1,Y_1,\cdots,Y_k)
\eeb
\beb
(X\wedge_H Y)(X_1,X_2,\cdots,X_k) = H(Y,X_1,X_3,\cdots X_{k})g(X,X_2)- H(X,X_1,X_3,\cdots X_{k})g(Y,X_2).
\eeb
%-------------------------------------------------
\indent A tensor $D\in \mathcal T^0_4(M)$ is said to be a generalized curvature tensor (\cite{DGHS11}, \cite{SDHJK15}, \cite{SK14}) if
$$D(X_1,X_2,X_3,X_4)+D(X_2,X_3,X_1,X_4)+D(X_3,X_1,X_2,X_4)=0,$$
$$D(X_1,X_2,X_3,X_4)+D(X_2,X_1,X_3,X_4)=0 \ \ \mbox{and}$$
$$D(X_1,X_2,X_3,X_4)=D(X_3,X_4,X_1,X_2).$$
We note that if $A,E \in \mathcal T^0_2(M)$ are both symmetric, then $A\wedge E$ is a generalized curvature tensor.
%------------------------------------
Again a generalized curvature tensor is called proper if 
$$(\nabla_{X_1} D)(X_2,X_3,X_4,X_5)+(\nabla_{X_2} D)(X_3,X_1,X_4,X_5)+(\nabla_{X_3} D)(X_1,X_2,X_4,X_5)=0.$$
%------------------------------------
Some important generalized curvature tensors are Gaussian curvature tensor $G$, conformal curvature tensor $C$, concircular curvature tensor $W$ and conharmonic curvature tensor $K$. These are respectively given as
\beb
G &=& \frac{1}{2} g \wedge g,\\
C &=& R -\frac{1}{n-2} g \wedge S+\frac{\kappa}{2(n-1)(n-2)}g \wedge g,\\
W &=& R - \frac{\kappa}{2 n(n-1)}g \wedge g \ \ \mbox{and}\\ 
K &=& R - \frac{1}{n-2} g \wedge S.
\eeb
%-------------------------------------
For $D\in \mathcal T^0_4(M)$ and $X,Y\in \chi(M)$, the associated $(1,3)$ tensor $\mathcal D$ and the associated curvature operator $\mathscr D(X,Y) \in \Xi(M)$ are respectively given by
$$g(\mathcal{D}(X,Y)X_1,X_2)=D(X,Y,X_1, X_2) \ \mbox{and}$$
$$\mathscr{D}(X,Y)(X_1)=\mathcal D(X,Y)X_1.$$
%------------------------------------------------
\indent One can operate an endomorphism $\mathscr L$ on a $(0,k)$ tensor $H$ and a $(1,k-1)$ tensor $\mathcal H$ (\cite{SK14}, \cite{SKgrt}) as
\beb
(\mathscr{L} H)(X_1,X_2,\cdots,X_k) &=& -H(\mathscr{L}X_1,X_2,\cdots,X_k) - \cdots - H(X_1,X_2,\cdots,\mathscr{L}X_k) \ \mbox{and}\\
(\mathscr{L} \mathcal H)(X_1,X_2,\cdots,X_{k-1}) &=& \mathscr{L}\mathcal H(X_1,X_2,\cdots,X_{k-1}) -\mathcal H(\mathscr{L}X_1,X_2,\cdots,X_{k-1})\\
&& - \cdots - \mathcal H(X_1,X_2,\cdots,\mathscr{L}X_{k-1})
\eeb
respectively. In particular, if we consider $\mathscr L = \mathscr D(X,Y)$ and $X\wedge_A Y$, then we get $D\cdot H, Q(A,H) \in \mathcal T^0_{k+2}(M)$ and $D\cdot \mathcal H, Q(A,\mathcal H) \in \mathcal T^1_{k+1}(M)$ respectively as follows:
\beb
&&D\cdot H(X_1,X_2, \ldots ,X_k,X,Y) = (\mathscr{D}(X,Y)\cdot H)(X_1,X_2, \ldots ,X_k)\\
&&= -H(\mathcal{D}(X,Y)X_1,X_2, \ldots ,X_k) - \cdots - H(X_1,X_2, \ldots ,\mathcal{D}(X,Y)X_k),
\eeb
\beb
&&Q(A,H)(X_1,X_2, \cdots ,X_k,X,Y) = ((X \wedge_A Y)\cdot H)(X_1,X_2, \ldots ,X_k)\\
&&= A(X, X_1) H(Y,X_2,\cdots,X_k) + \cdots + A(X, X_k) H(X_1,X_2,\cdots,Y)\\
&&- A(Y, X_1) H(X,X_2,\cdots,X_k) - \cdots - A(Y, X_k) H(X_1,X_2,\cdots,X),
\eeb
\beb
&&D\cdot \mathcal H(X_1,X_2, \cdots ,X_{k-1},X,Y) = (\mathscr{D}(X,Y)\cdot \mathcal H)(X_1,X_2, \cdots ,X_{k-1})\\
&&= \mathcal{D}(X,Y)\mathcal H(X_1,X_2, \cdots ,X_{k-1}) -\mathcal H(\mathcal{D}(X,Y)X_1,X_2, \cdots ,X_{k-1})\\
&&\hspace{0.4cm} - \cdots - \mathcal H(X_1,X_2, \cdots ,\mathcal{D}(X,Y)X_{k-1}),
\eeb
\beb
&&Q(A,\mathcal H)(X_1,X_2, \ldots ,X_{k-1},X,Y) = ((X \wedge_A Y)\cdot \mathcal H)(X_1,X_2, \cdots ,X_{k-1})\\
&&= A(Y,\mathcal H(X_1,X_2,\cdots,X_{k-1}))X-A(X,\mathcal H(X_1,X_2,\cdots,X_{k-1}))Y\\
&&+ A(X, X_1) \mathcal H(Y,X_2,\cdots,X_{k-1}) + \cdots + A(X, X_{k-1}) \mathcal H(X_1,X_2,\cdots,Y)\\
&&- A(Y, X_1) \mathcal H(X,X_2,\cdots,X_{k-1}) - \cdots - A(Y, X_{k-1}) \mathcal H(X_1,X_2,\cdots,X).
\eeb
%------------------------------------------------
\begin{defi} 
For $H \in \mathcal T^0_k(M)$ (resp., $\mathcal H \in \mathcal T^1_{k-1}(M)$) and $D \in \mathcal T^0_4(M)$, a semi-Riemannian manifold $M$ is said to be $H$-semisymmetric type (resp., $\mathcal H$-semisymmetric type) (\cite{SK14}, \cite{Szab82}) manifold due to $D$ if $D\cdot H = 0$ (resp., $D\cdot \mathcal H = 0$).
\end{defi}
In particular, a semi-Riemannian manifold respectively satisfying 
$R\cdot R = 0$, 
$R\cdot S = 0$, 
$R\cdot P = 0$, 
$P\cdot R = 0$ and 
$P\cdot S = 0$
is respectively called semisymmetric \cite{Szab82}, Ricci semisymmetric, projective semisymmetric, semisymmetric due to projective curvature tensor and Ricci semisymmetric due to projective curvature tensor.
%----------------------------------------------
\begin{defi}$($\cite{AD83}, \cite{Desz92}, \cite{DGHS11}, \cite{SK14}$)$ 
For $H \in \mathcal T^0_k(M)$ (resp., $\mathcal H \in \mathcal T^1_{k-1}(M)$) and $D_i \in \mathcal T^0_4(M)$, $i=1,2,\cdots r, r\ge 2$, a semi-Riemannian manifold is said to be $H$-pseudosymmetric type (resp., $\mathcal H$-pseudosymmetric type) if 
$\sum\limits_{i=1}^r c_i (D_i\cdot H) = 0$ \Big(resp., $\sum\limits_{i=1}^r c_i (D_i\cdot \mathcal H) = 0$\Big) for some $c_i\in C^{\infty}(M)$, called the associated scalars. Moreover a pseudosymmetric type condition is called constant type if its associated scalars are all constants.
\end{defi}
In particular, a semi-Riemannian manifold respectively satisfying
$R\cdot R = L_R Q(g,R)$, 
$R\cdot S = L_S Q(g,S)$, 
$R\cdot P = L_P Q(g,P)$, 
$P\cdot R = L_1 Q(g,R)$ and 
$P\cdot S = L_2 Q(g,S)$
is respectively called pseudosymmetric, Ricci pseudosymmetric, projective pseudosymmetric, pseudosymmetric due to projective curvature tensor and Ricci pseudosymmetric due to projective curvature tensor, where $L_R, L_S,L_P,L_1,L_2$ are the associated scalars.\\
%----------------------------------------------
\indent As a generalization of manifold of vanishing conformal curvature tensor (i.e., $C\equiv 0$ on $M$), there arose two curvature conditions, namely, Roter type \cite{Desz03} and generalized Roter type \cite{SKgrt}, which are respectively given by
\be\label{rt}
R=c_1 g\wedge g + c_2 g\wedge S + c_3 S\wedge S,
\ee
\be\label{grt}
R=c_1 g\wedge g + c_2 g\wedge S + c_3 S\wedge S + c_4 g\wedge S^2 + c_5 S\wedge S^2 + c_6 S^2\wedge S^2,
\ee
where $c_i \in C^{\infty}(M)$, $i=1,2,\cdots, 6$.
\begin{defi}
A semi-Riemannian manifold $M$ satisfying \eqref{rt} (resp., \eqref{grt}) for some $c_i \in C^{\infty}(M)$ is called a Roter type manifold $($cite{Desz03}, \cite{Desz03a}, \cite{DGHS11}, \cite{DGPV11}, \cite{DPS13} and \cite{Glog07}$)$ (resp., generalized Roter type manifold $($\cite{DGJPZ13}, \cite{DGJZ16}, \cite{DGPV15}, \cite{SDHJK15}, \cite{SKgrt}, \cite{SKgrtw} and \cite{Sawi15}$)$).
\end{defi}
\noindent\textbf{Note:} We note that every Roter type manifold is generalized Roter type and every manifold of vanishing conformal curvature tensor is Roter type. We also note that an Einstein generalized Roter type manifold is of constant curvature \cite{SKgrt}.
%-------------------------------------------------------------------
\begin{defi}$($\cite{Prav95}, \cite{Venz85}$)$ 
Let $\mathcal L(M)$ be the vector space formed by all 1-forms $\Theta$ on $M$ satisfying
$$\Theta(X_1)D(X_2,X_3,X_4,X_5)+\Theta(X_2)D(X_3,X_1,X_4,X_5)+\Theta(X_3)D(X_1,X_2,X_4,X_5) = 0,$$
where $D\in\mathcal T^0_4(M)$. Then $M$ is said to be a $D$-space by Venzi if $dim \mathcal L(M) \ge 1$.
\end{defi}
In \cite{Venz85} Venzi named such a space as $B$-space for $D=R$. 
%%%%%%%%%%%%%%%%%%%%%%%%%%%%%%%%%%%%%%%%%%%%%%%%%%%%%%%%%%%%%%%%%%%%%%%%%%%%%%%%%%%%%%%%%%%%%%%%%%%%%%%%%%%%%%%%%%%%
%                                               Basic results
%%%%%%%%%%%%%%%%%%%%%%%%%%%%%%%%%%%%%%%%%%%%%%%%%%%%%%%%%%%%%%%%%%%%%%%%%%%%%%%%%%%%%%%%%%%%%%%%%%%%%%%%%%%%%%%%%%%%
\section{Some curvature related properties}
%==========================================
\noindent From definitions we can state the following lemmas:
\begin{lem}\label{lem3.1}
If $H\in\mathcal T^0_k(M)$ and $A,E\in\mathcal T^0_2(M)$ are symmetric, then\\
(i) $X\wedge_H Y =0$ if and only if $H=0$,\\
(ii) $g \wedge H =0$ if and only if $H=0$ and\\
(iii) $Q(A,E)=0$ if and only if $A$ and $E$ are linearly dependent \cite{DD91}.
\end{lem}
\begin{lem}\label{lem3.2} 
If $A\in\mathcal T^0_2(M)$ is symmetric and $D\in\mathcal T^0_4(M)$ is a generalized curvature tensor, then 
$$D\cdot (X_1\wedge_A X_2) = X_1\wedge_{D\cdot A} X_2.$$ 
\end{lem}
\begin{lem}\label{lem3.3} \cite{SKgrt}
If $A,E\in\mathcal T^0_2(M)$ and $\mathscr L \in\Xi(M)$, then 
$$\mathscr L (A\wedge E) = A\wedge \mathscr L E + E\wedge \mathscr L A.$$ 
In particular, for $D\in\mathcal T^0_4(M)$ we have
$$D\cdot (A\wedge E) = A\wedge (D\cdot E) + E\wedge (D\cdot A).$$ 
\end{lem}
\begin{lem}\label{lem3.5}
If $A\in\mathcal T^0_2(M)$ is symmetric and $\Pi\in \chi^*(M)$, then
$$\Pi(X_1)(X_2\wedge_A X_3)(X,Y) +\Pi(X_2)(X_3\wedge_A X_1)(X,Y) + \Pi(X_3)(X_1\wedge_A X_2)(X,Y) = 0$$
holds if and only if \ $\Pi(X_1)A(X_2,X_3) = \Pi(X_2)A(X_1,X_3)$.
\end{lem}
\begin{lem}\label{lem3.4} \cite{DDV89}
If $A$ is a symmetric (0,2) tensor and $D$ is a generalized curvature tensor, then
\beb
(X_1\wedge_A X_2)X_3 \ + \ (X_2\wedge_A X_3)X_1 &+& (X_3\wedge_A X_1)X_2 = 0 \ \mbox{and}\\
Q(A,D)(X_1,X_2,X_3,X_4,X_5,X_6)&+&Q(A,D)(X_3,X_4,X_5,X_6,X_1,X_2)\\
&+&Q(A,D)(X_5,X_6,X_1,X_2,X_3,X_4)=0.
\eeb
\end{lem}
%
%=========================================
\begin{lem}\label{lem3.6}
On a semi-Riemannian manifold\\
(i) \bea\label{*1}
Q(g,P)(X_1,X_2,X_3,X_4,X_5,X_6)&+&Q(g,P)(X_3,X_4,X_5,X_6,X_1,X_2)\\\nonumber
&+&Q(g,P)(X_5,X_6,X_1,X_2,X_3,X_4)=0.
\eea
holds if and only if the manifold is Einstein.\\
(ii) \bea\label{*2}
Q(S,P)(X_1,X_2,X_3,X_4,X_5,X_6)&+&Q(S,P)(X_3,X_4,X_5,X_6,X_1,X_2)\\\nonumber
&+&Q(S,P)(X_5,X_6,X_1,X_2,X_3,X_4)=0.
\eea
holds if and only if $\kappa(n S - \kappa g) = 0$.
\end{lem}
\noindent \textbf{Proof:} (i) Contracting \eqref{*1} over $X_1$ and $X_4$, we get
\beb
&&\frac{\kappa g(X_5,X_3) g(X_6,X_2)}{n-1}-\frac{g(X_6,X_2)S(X_5,X_3)}{n-1}-\frac{g(X_2,X_3) S(X_5,X_6)}{n-1}\\
&&-\frac{n g(X_5,X_3) S(X_6,X_2)}{n-1}+\frac{g(X_5,X_3)S(X_6,X_2)}{n-1}+\frac{S(X_2,X_3) g(X_5,X_6)}{n-1} = 0.
\eeb
Again contracting the above over $X_2$ and $X_3$, we get
\beb
\frac{2 \kappa g(X_5,X_6)}{n-1}-\frac{2 n S(X_5,X_6)}{n-1} = 0 \ \Rightarrow \ n S(X_5,X_6) = \kappa g(X_5,X_6),
\eeb
which implies $M$ is Einstein.
(ii) Contracting \eqref{*2} over $X_1$ and $X_4$, we get
\bea\label{eq3.1}
&&\frac{\kappa g(X_2,X_6) S(X_3,X_5)}{n-1}-\frac{g(X_3,X_6)S^2(X_2,X_5)}{n-1}+\frac{g(X_3,X_5) S^2(X_2,X_6)}{n-1}\\\nonumber
&&-\frac{g(X_2,X_5)S^2(X_3,X_6)}{n-1}-\frac{g(X_2,X_3) S^2(X_5,X_6)}{n-1}-\frac{n S(X_2,X_6) S(X_3,X_5)}{n-1}\\\nonumber
&&-\frac{S(X_2,X_6) S(X_3,X_5)}{n-1}+\frac{2 S(X_2,X_5) S(X_3,X_6)}{n-1}+\frac{S(X_2,X_3) S(X_5,X_6)}{n-1} = 0.
\eea
Now contracting \eqref{eq3.1} over $X_2$ and $X_3$, we get
\beb
\frac{2 \kappa S(X_5,X_6)}{n-1}-\frac{2 n S^2(X_5,X_6)}{n-1} = 0 \ \Rightarrow \ S^2 = \frac{\kappa}{n} S \ \mbox{and} \ \kappa^{(2)} = \frac{\kappa^2}{n}.
\eeb
Again contracting \eqref{eq3.1} over $X_2$ and $X_5$, we get
\beb
-\frac{\kappa^{(2)} g(X_3,X_6)}{n-1}+\frac{3 \kappa S(X_3,X_6)}{n-1}-\frac{2 n S^2(X_3,X_6)}{n-1} = 0
\eeb
Thus putting the value of $S^2$ and $\kappa^{(2)}$ in the preceding equation, we get
$$\frac{\kappa (n S(X_3,X_6)-\kappa g(X_3,X_6))}{n(n-1)} = 0 \ \Rightarrow \ \kappa(n S(X_3,X_6) - \kappa g(X_3,X_6)) = 0.$$
%-----------------------------
\begin{pr}\label{prop3.1}
On a semi-Riemannian manifold, contraction and projective operator $\mathscr P(X,Y)$ commute if and only if the manifold is Einstein.
\end{pr}
\noindent \textbf{Proof:} The result follows from Lemma 5.1 and 5.2 of \cite{SK14}.
%-------------------------
\begin{cor}
$P\cdot G =0$ if and only if $M$ is Einstein.
\end{cor}
%--------------------------
\begin{pr}\label{prop3.2}
(i) $P\cdot S =R\cdot S$, \ \ (ii) $P\cdot (S\wedge S) =R\cdot (S\wedge S)$.
\end{pr}
%-------------------------
\begin{cor}
On a semi-Riemannian manifold, $P\cdot S = L Q(g,S) \Leftrightarrow R\cdot S = L Q(g,S)$, where $L$ is a smooth function on $\{x\in M:S_x\ne\frac{\kappa}{n}g_x\}$.
\end{cor}
%-------------------------
\begin{pr}\label{prop3.3}
On a semi-Riemannian manifold, the projective curvature tensor $P$ possesses the following identities:\\
(i) $\mathscr P(X,Y) = -\mathscr P(Y,X)$,\\
(ii) $\mathcal P(X_1,X_2)X_3 + \mathcal P(X_2,X_3)X_1 + \mathcal P(X_3,X_1)X_2 = 0$ and\\
(iii) $P(X_1,X_2,X_3,X) + P(X_2,X_3,X_1,X) + P(X_3,X_1,X_2,X) = 0$.
\end{pr}
%-------------------------
\begin{pr}\label{prop3.4}
On a semi-Riemannian manifold the following conditions are equivalent:\\
(i) $P(X_1,X_2,X_3,X_4) + P(X_1,X_2,X_4,X_3) = 0$,\\
(ii) $P(X,X_1,X_2,X_3) + P(X,X_2,X_3,X_1) + P(X,X_3,X_1,X_2) = 0$ and\\
(iii) $M$ is Einstein.
\end{pr}
%-------------------------
\begin{pr}\label{prop3.5}
On a semi-Riemannian manifold $M$,\\
(i) $(\nabla_{X_1} P)(X_2,X_3,X,Y) + (\nabla_{X_2} P)(X_3,X_1,X,Y) + (\nabla_{X_3} P)(X_1,X_2,X,Y) = 0$ if and only if\\ 
$(\nabla_{X_1}S)(X_2,X_3)-(\nabla_{X_2}S)(X_1,X_3) = 0$, i.e., the Ricci tensor of $M$ is Codazzi type and\\
(ii) $(\nabla_{X_1} P)(X,Y,X_2,X_3) + (\nabla_{X_2} P)(X,Y,X_3,X_1) + (\nabla_{X_3} P)(X,Y,X_1,X_2) = 0$ if and only if\\ 
$(\nabla_{X_1}S)(X_2,X_3) = 0$, i.e., $M$ is Ricci symmetric.
\end{pr}
\noindent \textbf{Proof:} The result (i) follows from Proposition 2.2 of \cite{SH09}, and (ii) can be proved in a similar way.
%-------------------------
\begin{pr}\label{prop3.6} \cite{SK14}
On a semi-Riemannian manifold, we have the following:\\
(i) $\nabla R =0$ $\Leftrightarrow$ $\nabla P = 0$ (see also \cite{RT67}), \\
(ii) $R\cdot R =0$ $\Leftrightarrow$ $R\cdot P = 0$, \\
(iii) $R\cdot R = L Q(g,R)$ $\Leftrightarrow$ $R\cdot P = L Q(g,P)$, \\
(iv) $D\cdot R =0$ $\Leftrightarrow$ $D\cdot P = 0$, \\
(v) $D\cdot R = L Q(g,R)$ $\Leftrightarrow$ $D\cdot P = L Q(g,P),$\\
where $D$ is a generalized curvature tensor and $L\in C^{\infty}(M)$.
\end{pr}
%--------------------------
\begin{pr}\label{prop3.7}
If $H\in \mathcal T^0_k(M)$ such that $Q(S,H) =0$, then\\
(i) $P\cdot H =0$ $\Leftrightarrow$ $R\cdot H =0$ and\\
(ii) $P\cdot H =L Q(g,H)$ $\Leftrightarrow$ $R\cdot H = L Q(g,H)$, $L\in C^{\infty}(M)$.
\end{pr}
%===============================================================
\begin{lem}
On a semi-Riemannian manifold the following conditions are equivalent:\\
(i) $n^2 S^2 - 2 n \kappa S + \kappa^2 g = 0$ \ \ \ (ii) $n S^2 - 2 \kappa S + \kappa^{(2)} g = 0$, \ where $\kappa^{(2)} = \mbox{trace}(\mathcal S^2)$.
\end{lem}
\noindent \textbf{Proof:} Since under contraction both the conditions give $\kappa^{(2)} = \frac{1}{n}\kappa^2$, hence the result is obvious.
%----------------------------------
\begin{pr}\label{pr4.6}
A Riemannian manifold satisfying the curvature condition $n^2 S^2 - 2 n \kappa S + \kappa^2 g = 0$ or $n S^2 - 2 \kappa S + \kappa^{(2)} g = 0$ is always an Einstein manifold.
\end{pr}
\noindent \textbf{Proof:} If $M$ is a Riemannian manifold, then at each point $x\in M$ the Ricci operator $\mathcal S$ is symmetric and there exists an orthonormal basis $\{e_1, e_2, \cdots, e_n\}$ of $(T_x M, g_x)$ consisting of eigenvectors of $\mathcal S_x$. Let $\mathcal S_x e_i = \lambda_i e_i$ for each $i = 1, 2, \cdots, n$, where $\lambda_i \in \mathbb R$ are the corresponding eigenvalues.\\
Therefore $S(e_i,e_i) = \lambda_i$ and $S^2(e_i,e_i) = \lambda_i^2$ for each $i$.
Now from the given curvature condition, we have
$$n^2 \lambda_i^2 - 2 n \kappa \lambda_i + \kappa^2 = 0, \forall \ i$$
$$\Rightarrow (n \lambda_i -\kappa)^2 = 0 \Rightarrow \lambda_i = \frac{\kappa}{n}, \forall \ i.$$
Hence all the eigenvalues are equal and thus $M$ is an Einstein manifold.
%==============================================================================
%%%%%%%%%%%%%%%%%%%%%%%%%%%%%%%%%%%%%%%%%%%%%%%%%%%%%%%%%%%%%%%%%%%%%%%%%%%%%%%%%%%%%%%%%%%%%%%%%%%%%%%%%%%%%%%%%%%%%%%%%
%%%%%%%%%%%%%%%%%%%%%%%%%%%%%%%%%%%%%%%%%%%%%%%%%%%%%%%%%%%%%%%%%%%%%%%%%%%%%%%%%%%%%%%%%%%%%%%%%%%%%%%%%%%%%%%%%%%%%%%%%
%                                                  Main results
%%%%%%%%%%%%%%%%%%%%%%%%%%%%%%%%%%%%%%%%%%%%%%%%%%%%%%%%%%%%%%%%%%%%%%%%%%%%%%%%%%%%%%%%%%%%%%%%%%%%%%%%%%%%%%%%%%%%%%%%%
\section{ Some pseudosymmetric type curvature conditions}
%%%%%%%%%%%%%%%%%%%%%%%%%%%%%%%%%%%%%%%%%%%%%%%%%%%%%%%%%%
%%%%%%%%%%%%%%%%%%%%%%%%%%%%%%%%%%%%%%%%%%%%%%%%%%%%%%%%%%%%%%%%%%%%%%%%
%                  Walker type conditions
%%%%%%%%%%%%%%%%%%%%%%%%%%%%%%%%%%%%%%%%%%%%%%%%%%%%%%%%%%%%%%%%%%%%%%%
It is well known that every semi-Riemannian manifold $M$ satisfies
\beb
R\cdot R(X_1,X_2,X_3,X_4,X_5,X_6) &+& R\cdot R(X_3,X_4,X_5,X_6,X_1,X_2)\\
&+& R\cdot R(X_5,X_6,X_1,X_2,X_3,X_4)=0.
\eeb
This identity is known as Walker identity. For two (0,4) tensors $D_1$ and $D_2$ on a semi-Riemannian manifold, the condition
\bea\label{wlc}
D_1\cdot D_2(X_1,X_2,X_3,X_4,X_5,X_6) &+& D_1\cdot D_2(X_3,X_4,X_5,X_6,X_1,X_2)\\\nonumber
&+& D_1\cdot D_2(X_5,X_6,X_1,X_2,X_3,X_4)=0
\eea
is called Walker type condition. If for some particular $D_1$ and $D_2$, the condition \eqref{wlc} holds identically on every semi-Riemannian manifold, then it is called an Walker type identity (\cite{DHS12}, \cite{DY07}).
%=========================================
\begin{pr}\label{pr4.1}
Every semi-Riemannian manifold satisfies the following Walker type identity:
\beb
D\cdot D(X_1,X_2,X_3,X_4,X_5,X_6) &+& D\cdot D(X_3,X_4,X_5,X_6,X_1,X_2)\\
&+& D\cdot D(X_5,X_6,X_1,X_2,X_3,X_4)=0,
\eeb
where $D$ is a generalized curvature tensor.
\end{pr}
%=========================================
Since $R\cdot W = R\cdot R$, $W\cdot R = R\cdot R -\frac{\kappa}{2n(n-1)}Q(g,R)$, $P\cdot R = R\cdot R -\frac{1}{n-1}Q(S,R)$, $K\cdot C = K\cdot K$ and $C\cdot K = K\cdot K -\frac{\kappa}{2 (n-1)(n-2)}Q(g,K)$  then in view of Proposition \ref{pr4.1} and Lemma \ref{lem3.4} we can state the following:
\begin{pr}\label{pr4.2}
Every semi-Riemannian manifold satisfies the following Walker type identities:
\beb
R\cdot W(X_1,X_2,X_3,X_4,X_5,X_6) &+& R\cdot W(X_3,X_4,X_5,X_6,X_1,X_2)\\
&+& R\cdot W(X_5,X_6,X_1,X_2,X_3,X_4)=0,
\eeb
\beb
W\cdot R(X_1,X_2,X_3,X_4,X_5,X_6) &+& W\cdot R(X_3,X_4,X_5,X_6,X_1,X_2)\\
&+& W\cdot R(X_5,X_6,X_1,X_2,X_3,X_4)=0,
\eeb
\beb
P\cdot R(X_1,X_2,X_3,X_4,X_5,X_6) &+& P\cdot R(X_3,X_4,X_5,X_6,X_1,X_2)\\
&+& P\cdot R(X_5,X_6,X_1,X_2,X_3,X_4)=0,
\eeb
\beb
C\cdot K(X_1,X_2,X_3,X_4,X_5,X_6) &+& C\cdot K(X_3,X_4,X_5,X_6,X_1,X_2)\\
&+& C\cdot K(X_5,X_6,X_1,X_2,X_3,X_4)=0,
\eeb
\beb
K\cdot C(X_1,X_2,X_3,X_4,X_5,X_6) &+& K\cdot C(X_3,X_4,X_5,X_6,X_1,X_2)\\
&+& K\cdot C(X_5,X_6,X_1,X_2,X_3,X_4)=0.
\eeb
\end{pr}
%%%%%%%%%%%%%%%%%%%%%%%%%%%%%%%%%%%%%%%%%%%%%%%%%%%%%%%%%%%%%%%%%%%%%%%%%%%%%%%%%%%%%%%%%%%%%%%%%%%%%%%%%%%%%%%%%%%%%%%%%
%===============================================================================
\begin{thm}\label{thm4.1}
On a semi-Riemannian manifold $M$,
\bea\label{r.p}
R\cdot P(X_1,X_2,X_3,X_4,X_5,X_6) &+& R\cdot P(X_3,X_4,X_5,X_6,X_1,X_2)\\\nonumber
&+& R\cdot P(X_5,X_6,X_1,X_2,X_3,X_4)=0
\eea
holds if and only if $M$ is Ricci semisymmetric.
\end{thm}
%-----------------------
\noindent \textbf{Proof:} Let us first consider $R\cdot S =0$. Then from Lemma \ref{lem3.2}, $R\cdot P = R\cdot R$ and hence \eqref{r.p} reduces to the Walker identity.\\
For the converse part, contracting \eqref{r.p} over $X_1$ and $X_3$, we get
\bea\label{eq4.3}
&&-R(X_2, X_6, X_4, \mathcal S(X_5)) - R(X_2, X_6, X_5, \mathcal S(X_4))\\\nonumber
&&+ R(X_2, \mathcal S(X_4), X_5, X_6) +  R(X_4, \mathcal S(X_2), X_5, X_6) \\\nonumber
&&+ g(X_4, X_6) \left[-E(X_2, X_5) + \mathcal S^2(X_2, X_5) \right]\\\nonumber
&&+ g(X_2, X_6) \left[-E(X_4, X_5) + \mathcal S^2(X_4, X_5)\right]\\\nonumber
&&+ g(X_2, X_5) \left[E(X_4, X_6) - \mathcal S^2(X_4, X_6)\right] = 0,
\eea
where $E(X,Y)$ is the tensor obtained from $R(X_1, X, Y, \mathcal S(X_2))$ by taking contraction over $X_1$ and $X_2$. Again contracting \eqref{r.p} over $X_1$ and $X_4$ and then replacing $X_3$ by $X_4$, we get
\bea\label{eq4.4}
&& R(X_2, X_6, X_4, \mathcal S(X_5)) +  R(X_2, X_6, X_5, \mathcal S(X_4)) \\\nonumber
&&+ (-1 + n) \left[R(X_2, \mathcal S(X_4), X_5, X_6) + R(X_4, \mathcal S(X_2), X_5, X_6)\right] \\\nonumber
&&+ g(X_4, X_6) \left[E(X_2, X_5) - \mathcal S^2(X_2, X_5)\right] \\\nonumber
&&+ g(X_2, X_6) \left[E(X_4, X_5) - \mathcal S^2(X_4, X_5)\right] \\\nonumber
&&+ g(X_2, X_5) \left[-E(X_4, X_6) + \mathcal S^2(X_4, X_6)\right] = 0.
\eea
Now adding \eqref{eq4.3} and \eqref{eq4.4}, we get $n\left[R(X_2, \mathcal S(X_4), X_5, X_6) + R(X_4, \mathcal S(X_2), X_5, X_6)\right]=0$, i.e., $R\cdot S = 0$. This completes the proof.
%=================================================================================
%==================================================================================
\begin{thm}\label{thm4.2}
On a semi-Riemannian manifold $M$,
\bea\label{p.pc}
P\cdot P(X_1,X_2,X_3,X_4,X_5,X_6) &+& P\cdot P(X_3,X_4,X_5,X_6,X_1,X_2)\\\nonumber
&+& P\cdot P(X_5,X_6,X_1,X_2,X_3,X_4)=0
\eea
holds if $\wedge_{R\cdot S} = \frac{1}{n-1} \wedge_S\cdot\wedge_S$. Again if \eqref{p.pc} holds then\\
(i) $n(n-1) R\cdot S = \kappa Q(g, S)$,\\
(ii) $S\wedge S = g\wedge S^2$,\\
(iii) $n^2 S^2 - 2 n \kappa S + \kappa^2 g = 0$. Moreover if the manifold is Riemannian, then it becomes an Einstein manifold.
\end{thm}
%------------------------
\noindent \textbf{Proof:} Let us first consider $\wedge_{R\cdot S} = \frac{1}{n-1} \wedge_S\cdot\wedge_S$. Then by Lemma \ref{lem3.2}, $R\cdot \wedge_{S}= \frac{1}{n-1} \wedge_S\cdot\wedge_S$. Thus in view of Proposition \ref{pr4.2}, the left hand side of \eqref{p.pc} reduces to
\beb
(P\cdot \wedge_S)(X_1,X_2,X_3,X_4,X_5,X_6) &+& (P\cdot \wedge_S)(X_3,X_4,X_5,X_6,X_1,X_2)\\
&+& (P\cdot \wedge_S)(X_5,X_6,X_1,X_2,X_3,X_4)
\eeb
\beb
&=& (R\cdot \wedge_S)(X_1,X_2,X_3,X_4,X_5,X_6) - \frac{1}{n-1} (\wedge_S\cdot\wedge_S)(X_1,X_2,X_3,X_4,X_5,X_6)\\
&&+ (R\cdot \wedge_S)(X_3,X_4,X_5,X_6,X_1,X_2) - \frac{1}{n-1} (\wedge_S\cdot\wedge_S)(X_3,X_4,X_5,X_6,X_1,X_2)\\
&&+ (R\cdot \wedge_S)(X_5,X_6,X_1,X_2,X_3,X_4) - \frac{1}{n-1} (\wedge_S\cdot\wedge_S)(X_5,X_6,X_1,X_2,X_3,X_4)\\
&=& 0
\eeb
Now similar to the proof of the converse part of Theorem \ref{thm4.1}, we get 
\be\label{eq4.6}
n(n-1) R\cdot S = \kappa Q(g, S) + (S\wedge S - g\wedge S^2).
\ee
Again contracting \eqref{p.pc} over $X_5$ and $X_6$, we get
\bea\label{eq4.7}
&&-\frac{g(X_2,X_4) S^2(X_1,X_3)}{(n-1)^2}+\frac{g(X_2,X_3) S^2(X_1,X_4)}{(n-1)^2}+\frac{g(X_1,X_4)S^2(X_2,X_3)}{(n-1)^2}\\\nonumber
&&-\frac{g(X_1,X_3) S^2(X_2,X_4)}{(n-1)^2}-\frac{2 S(X_1,X_4) S(X_2,X_3)}{(n-1)^2}+\frac{2 S(X_1,X_3) S(X_2,X_4)}{(n-1)^2} = 0.
\eea
Rearranging we get $\frac{1}{(n-1)^2}\left[(g\wedge S^2)(X_1,X_2,X_3,X_4)-(S\wedge S)(X_1,X_2,X_3,X_4)\right]$. Hence from \eqref{eq4.6}, we get $n(n-1) R\cdot S = \kappa Q(g, S)$.\\
Again contracting \eqref{eq4.7}, we get
$$\frac{1}{(n-1)^2}\left[\kappa^{(2)} g(X_2,X_3)+n S^2(X_2,X_3)-2 \kappa S(X_2,X_3)\right] = 0.$$
Then the theorem directly follows from Proposition \ref{pr4.6}.
%========================================================================
\begin{rem}
Although the Walker type condition presented in Theorem \ref{thm4.1} is equivalent to Ricci semisymmetry but the Walker type condition presented in Theorem \ref{thm4.2} implies $R\cdot R = \frac{\kappa}{n(n-1)} Q(g,S)$ but not conversely. To support this fact in Example 3 we present a 5-dimensional semi-Riemannian manifold which satisfies $R\cdot R = \frac{\kappa}{n(n-1)} Q(g,S)$ and $$P\cdot P(X_1,X_2,X_3,X_4,X_5,X_6) + P\cdot P(X_3,X_4,X_5,X_6,X_1,X_2) + P\cdot P(X_5,X_6,X_1,X_2,X_3,X_4) \neq 0.$$
\end{rem}
%%%%%%%%%%%%%%%%%%%%%%%%%%%%%%%%%%%%%%%%%%%%%%%%%%%%%%%%%%%%%%%%%%%%%%%%%%%%%%%%%%%%%%%%%%%%%%%%%%%%%%%%%%%%%%%%%%%%%%%%%
%                                                  P-space by Venzi
%%%%%%%%%%%%%%%%%%%%%%%%%%%%%%%%%%%%%%%%%%%%%%%%%%%%%%%%%%%%%%%%%%%%%%%%%%%%%%%%%%%%%%%%%%%%%%%%%%%%%%%%%%%%%%%%%%%%%%%%%
\begin{pr}\label{pr4.3}
Let $M$ be a semi-Riemannian $P$-space by Venzi with associated 1-form $\Pi$. Then it is a\\
(i) $R$-space by Venzi if and only if 
$$\Pi(X_2)S(X_1,X_3) - \Pi(X_1)S(X_2,X_3)=0,$$
(ii) $W$-space by Venzi if and only if 
$$\Pi(X_2)\left(S-\frac{\kappa}{n}g\right)(X_1,X_3) - \Pi(X_1)\left(S-\frac{\kappa}{n}g\right)(X_2,X_3)=0.$$
Moreover in both cases $\kappa = 0$.
\end{pr}
\noindent \textbf{Proof:} Since $R-P=\frac{1}{n-1}\wedge_S$ and $W-P=\frac{1}{n(n-1)}(\wedge_{n S-\kappa g})$, the results follows from Lemma \ref{lem3.5}.

%=================================================================================
\begin{thm}\label{thm4.3}
Let $M$ be a semi-Riemannian $P$-space by Venzi with associated 1-form $\Pi$. If\\
(i) $\Pi$ is non-null at $x\in M$, then $W=0$ at $x$.\\
(ii) $\Pi$ is null at some $x\in M$, then $M$ is a $W$-space by Venzi with same associated 1-form.
\end{thm}
%--------------------
\noindent \textbf{Proof:} From hypothesis
\be\label{eq4.9}
\Pi(X_3)P(X_1,X_2,X_4,X_5)+\Pi(X_2)P(X_3,X_1,X_4,X_5)+\Pi(X_1)P(X_2,X_3,X_4,X_5)=0.
\ee
Let $V$ be the vector field corresponding to $\Pi$. Then contracting \eqref{eq4.9} over $X_3$ and $X_5$, we get
$$P(X_1,X_2,X_4,V) = 0.$$
\beb
\Rightarrow R(X_1,X_2,X_4,V) &=& \frac{1}{n-1}\left[S(X_2,X_4)g(X_1,V)-S(X_1,X_4)g(X_2,V)\right]\\
											&=& \frac{1}{n-1}\left[S(X_2,X_4)\Pi(X_1)-S(X_1,X_4)\Pi(X_2)\right]\\
\Rightarrow S(V,X_1) &=& \frac{1}{n-1}\left[\kappa \Pi(X_1)-S(X_1,V)\right].
\eeb
Hence
\beb
R(V,X_4,X_1,X_2) &=& \frac{1}{n-1}\left[S(X_1,X_4)\Pi(X_2) - S(X_2,X_4)\Pi(X_1)\right],\\
S(V,X_1) &=& \frac{\kappa}{n}\Pi(X_1) \ \ \mbox{and}
\eeb
\bea\label{Pv}
P(V,X_4,X_1,X_2) &=& \frac{1}{n-1}\left[S(X_1,X_4)\Pi(X_2) - S(X_2,X_4)\Pi(X_1)\right]\\\nonumber
									& & - \frac{1}{n-1}\left[S(X_1,X_4)\Pi(X_2) - S(V,X_1)g(X_2,X_4)\right]\\\nonumber
								&=& \frac{1}{n-1}\left[\frac{\kappa}{n}g(X_2,X_4) - S(X_2,X_4)\right]\Pi(X_1).
\eea
%---------------------
(i) Let us suppose that $\Pi$ is non-null at $x$, and without loss of generality, we can consider the associated vector field $V$ of $\Pi$ is of unit norm. 
Now putting $X_3 = V$ in \eqref{eq4.9}, we get
\beb
P(X_1,X_2,X_4,X_5) &=& -\Pi(X_2)P(V,X_1,X_4,X_5) + \Pi(X_1)P(V,X_2,X_4,X_5)\\
									&=& \frac{\Pi(X_2)\Pi(X_4)}{n-1}\left(S-\frac{\kappa}{n}g\right)(X_1,X_5) - \frac{\Pi(X_1)\Pi(X_4)}{n-1}\left(S-\frac{\kappa}{n}g\right)(X_2,X_5)\\
									&=& \frac{\Pi(X_4)}{n-1}\left[\Pi(X_2)\left(S-\frac{\kappa}{n}g\right)(X_1,X_5) - \Pi(X_1)\left(S-\frac{\kappa}{n}g\right)(X_2,X_5)\right]\\
									&=& \frac{1}{n-1}\Pi(X_4)\left[\Pi(X_2)Z(X_1,X_5) - \Pi(X_1)Z(X_2,X_5)\right],
\eeb
where $Z=S-\frac{\kappa}{n}g$. Thus the curvature tensor $R$ is given by
\bea\label{eq4.8}
R(X_1,X_2,X_4,X_5) &=& \frac{1}{n-1}\Pi(X_4)\left[\Pi(X_2)Z(X_1,X_5) - \Pi(X_1)Z(X_2,X_5)\right] \\\nonumber
										& & + \frac{1}{n-1}\left[S(X_2,X_4)g(X_1,X_5) - S(X_1,X_4)g(X_2,X_5)\right].
\eea
As $Z$ is trace free and $Z(X,V) = 0$, $\forall\ X$, contracting \eqref{eq4.8} over $X_1$ and $X_5$, we get
$$S(X_2,X_4) = -\frac{1}{n-1}(\kappa g(X_2,X_4) + S(X_2,X_4))$$
$$\Rightarrow S(X_2,X_4) - \frac{\kappa}{n}g(X_2,X_4) = 0 \ \Rightarrow \ Z(X_2,X_4) =0.$$
Again putting this in \eqref{eq4.8}, we get $R(X_1,X_2,X_4,X_5) = \frac{\kappa}{n(n-1)}G(X_1,X_2,X_4,X_5)$. This completes the proof of (i).\\
%=============================
(ii) If $\Pi$ is null at $x$, then $\Pi(V) = 0$. Putting $X_3 = V$ in \eqref{eq4.9} and using \eqref{Pv}, we get 
$$\Pi(X_2)\left(S-\frac{\kappa}{n}g\right)(X_1,X_5) - \Pi(X_1)\left(S-\frac{\kappa}{n}g\right)(X_2,X_5)=0.$$
Hence from Proposition \ref{pr4.3}, we get
$$\Pi(X_3)W(X_1,X_2,X_4,X_5)+\Pi(X_2)W(X_3,X_1,X_4,X_5)+\Pi(X_1)W(X_2,X_3,X_4,X_5)=0$$
at $x$. Now by (i), if $\Pi$ is non-null at $x$ then $W=0$ at $x$ and the above condition is obvious. Hence $M$ is a $W$-space by Venzi with same associated 1-form $\Pi$. This completes the proof.
%================================================================
%
\begin{pr}\label{pr4.3a} (Theorem 1, \cite{DG90})
If a generalized curvature tensor $D$ satisfies
$$\Pi(X_3)D(X_1,X_2,X_4,X_5)+\Pi(X_2)D(X_3,X_1,X_4,X_5)+\Pi(X_1)D(X_2,X_3,X_4,X_5)=0$$
for a 1-form $\Pi$, then $D\cdot D = Q(Ric(D), D)$ at the points where $\Pi \neq 0$ ($Ric(D)$ is the trace of the linear map $X_1\rightarrow \mathcal D(X_1,X_2)X_3$, i.e., the Ricci tensor corresponding to $D$).
\end{pr}

\indent Now in view of Theorem \ref{thm4.3} and Proposition \ref{pr4.3a} we can state the following:
\begin{thm}\label{thm4.4}
Every semi-Riemannian $P$-space by Venzi satisfies $W\cdot W = Q(S-\frac{\kappa}{n}g, W)$.
\end{thm}
Since on a Riemannian manifold all the non-zero 1-forms are non-null and on a manifold of constant curvature, $P =0$, hence we can state the following:
\begin{cor}\label{cor4.1}
Every Riemannian manifold is a $P$-space by Venzi if and only if it is of constant curvature.
\end{cor}

\begin{rem}
From Theorem \ref{thm4.3}, \ref{thm4.4} and Corollary \ref{cor4.1} it is clear that there does not exist any proper $P$-space by Venzi Riemannian manifold but such a structure exists on a semi-Riemannian manifold with a null associated 1-form. To support this fact in Example 3 we present a 5-dimensional semi-Riemannian manifold which is a $P$-space by Venzi.
\end{rem}
%%%%%%%%%%%%%%%%%%%%%%%%%%%%%%%%%%%%%%%%%%%%%%%%%%%%%%%%%%%%%%%%%%%%%%%%%%%%%%%%%%%%%%%%%%%%%%%%%%%%%%%%%%%%%%%%%%%%%%%%%
%                                P.R(0,4)=P.R(1,3), P.S(0,2)=P.S(1,1), P.P(0,4)=P.P(1,3)
%%%%%%%%%%%%%%%%%%%%%%%%%%%%%%%%%%%%%%%%%%%%%%%%%%%%%%%%%%%%%%%%%%%%%%%%%%%%%%%%%%%%%%%%%%%%%%%%%%%%%%%%%%%%%%%%%%%%%%%%%
\indent It is well known that if $D$ is a generalized curvature tensor, then 
\beb
g((D\cdot\mathcal R)(X_1,X_2,X_3,X,Y),X_4) &=& (D\cdot R)(X_1,X_2,X_3,X,Y,X_4) \ \mbox{and}\\ 
g((D\cdot\mathcal S)(X_1,X,Y),X_2) &=& (D\cdot S)(X_1,X_2,X,Y)
\eeb
But these results are not true for $D=P$. In this case we have the following:
%==============================================================================================
\begin{thm}\label{thm4.5}
A semi-Riemannian manifold $M$ satisfies 
\be\label{eqp.s}
g(P\cdot\mathcal S(X_1,X,Y),X_2) = P\cdot S(X_1,X_2,X,Y)
\ee
if and only if $(S\wedge S)(X_1, X_2, X,Y) = 2 (X\wedge_{S^2}Y)(X_1, X_2)$. Moreover if $M$ satisfies \eqref{eqp.s}, then $n S^2 = \kappa S = \frac{\kappa^2}{n} g$ and hence $M$ is Einstein if $\kappa$ is non-zero on $\{x\in M: \left(S-\frac{\kappa}{n} g\right)_x \neq 0\}$.
\end{thm}
%------------------------
\noindent \textbf{Proof:} From definition, we have
$$
P\cdot S(X_1,X_2,X,Y) = -S(\mathcal P(X,Y)X_1, X_2) - S(X_1, \mathcal P(X,Y)X_2)
$$
$$
\mbox{and } P\cdot \mathcal S(X_1,X,Y) = \mathcal P(X,Y)\mathcal S(X_1) - \mathcal S(\mathcal P(X,Y)X_1).
$$
Hence $g((P\cdot\mathcal S)(X_1,X,Y),X_2) = P(X,Y,\mathcal S(X_1),X_2)- S(\mathcal P(X,Y)X_1, X_2)$.\\
So $g(P\cdot\mathcal S(X_1,X,Y),X_2) = P\cdot S(X_1,X_2,X,Y)$ holds if and only if
$$P(X,Y,\mathcal S(X_1),X_4) = -S(X_1,\mathcal P(X,Y)X_2).$$
\be\label{s.s}
\Leftrightarrow (S\wedge S)(X_1, X_2, X,Y) = 2 (X\wedge_{S^2}Y)(X_1, X_2).
\ee
Now contracting \eqref{s.s} over $X$ and $X_2$, we get $n S^2 = \kappa S$. Again contracting \eqref{s.s} over $Y$ and $X_1$, we get $\kappa S = \kappa^{(2)}g$, which implies $\kappa^{(2)} = \frac{\kappa^2}{n}$. This completes the proof.
%=============================================================================================
%
\begin{rem}
From the above theorem we can say that if a manifold satisfies \eqref{eqp.s} with $\kappa = 0$ then the manifold may or may not be Einstein. To support this result, we present a non-Einstein semi-Riemannian metric with zero scalar curvature in Example 4, which satisfies $P\cdot S  = 0$ as well as $P\cdot \mathcal S = 0$.
\end{rem}
%
%=============================================================================================
\begin{thm}\label{thm4.6}
A semi-Riemannian manifold $M$ satisfies
\be\label{eqp.r}
g((P\cdot\mathcal R)(X_1,X_2,X_3,X,Y),X_4)=P\cdot R(X_1,X_2,X_3,X_4,X,Y),
\ee
if and only if $M$ is Einstein.
\end{thm}
%---------------------------------------------------------
\noindent \textbf{Proof:} From \eqref{eqp.r}, we get
\bea\label{eqp.ri}
& & g(X_4, Y) R(X_1, X_2, X_3, \mathcal S X) - g(X, X_4) R(X_1, X_2, X_3, \mathcal S Y) \\\nonumber
&=& - R(X_1, X_2, X_3, Y) S(X, X_4) - R(X, X_3, X_1, X_2) S(X_4, Y).
\eea
Now contracting \eqref{eqp.ri} over $X$ and $X_4$, we get $R(X_1,X_2,X_3,\mathcal S Y) = \frac{\kappa}{n} R(X_1,X_2,X_3,Y)$. Now putting the value of $R(X_1,X_2,X_3,\mathcal S Y)$ in \eqref{eqp.ri}, we get
\beb
& & \kappa [g(X_4, Y) R(X_1, X_2, X_3, X) - g(X, X_4) R(X_1, X_2, X_3, Y)] \\\nonumber
&=& - n[R(X_1, X_2, X_3, Y) S(X, X_4) + R(X, X_3, X_1, X_2) S(X_4, Y)].
\eeb
Now contracting the above equation over $X$ and $X_1$ and putting the value of $R(X_1,X_2,X_3,\mathcal S Y)$, we get
$$S(X_2,X_3)[n S(X_4,Y) - \kappa g(X_4,Y)] = 0.$$
This implies that the manifold is Einstein. We know that on an Einstein manifold, $P = W$, hence the converse part is obvious.
%==============================================================================================
%
%
%=======================================================================================
\begin{thm}\label{thm4.7}
A semi-Riemannian manifold satisfies
\be\label{eqp.p}
g((P\cdot\mathcal P)(X_1,X_2,X_3,X,Y),X_4)=P\cdot P(X_1,X_2,X_3,X_4,X,Y),
\ee
if and only if $M$ is Einstein.
\end{thm}
%--------------------------------------------------------------------
\noindent \textbf{Proof:} Contacting \eqref{eqp.p} over $X$ and $X_4$, we get
\bea\label{eqp.pi}
&& R(X_1,X_2,X_3,\mathcal S Y) = \frac{\kappa}{n} R(X_1,X_2,X_3,Y)\\\nonumber
                       &&\indent\hspace{1in}   + \frac{1}{2(n-1)}(S\wedge S)(X_1,X_2,X_3,Y)-\frac{\kappa}{n(n-1)}g((X_1\wedge_S X_2)X_3,Y).
\eea
Now contracting \eqref{eqp.p} over $X_1$ and $Y$ and using \eqref{eqp.pi}, we get
\bea\label{eqp.pi1}
-\frac{\kappa g(X_3,X_4) S(X,X_2)}{(n-1)^2 n}&+&\frac{\kappa g(X,X_2) S(X_3,X_4)}{(n-1)^2 n}\\\nonumber
&=&\frac{g(X,X_2) S^2(X_3,X_4)}{(n-1)^2}-\frac{S(X_3,X_4) S(X,X_2)}{(n-1)^2}
\eea
\be\label{s2}
\Rightarrow n^2 S^2(X_3, X_4) = 2 n \kappa S(X_3, X_4) - \kappa^2 g(X_3, X_4).
\ee
Again contracting \eqref{eqp.p} over $X_3$ and $Y$ and using \eqref{eqp.pi}, we get
\bea\label{eqp.pi2}
&& \frac{\kappa g(X_2,X_4) S(X,X_1)}{n (n-1)^2} - \frac{\kappa g(X_1,X_4) S(X,X_2)}{n (n-1)^2} + \frac{\kappa g(X,X_2) S(X_1,X_4)}{n (n-1)^2}\\\nonumber
&& -\frac{\kappa g(X,X_1) S(X_2,X_4)}{n (n-1)^2}-\frac{g(X,X_2) S^2(X_1,X_4)}{(n-1)^2}+\frac{g(X,X_1) S^2(X_2,X_4)}{(n-1)^2}\\\nonumber
&& -\frac{S(X_2,X_4) S(X,X_1)}{(n-1)^2}+\frac{S(X,X_2) S(X_1,X_4)}{(n-1)^2} = 0.
\eea
Now replacing $X_1$ by $X_3$ in \eqref{eqp.pi2} and then subtracting from \eqref{eqp.pi1}, we get
$$\frac{\kappa g(X,X_3) S(X_2,X_4)}{n(n-1)^2}+\frac{S(X_2,X_4) S(X,X_3)}{(n-1)^2} = \frac{\kappa g(X_2,X_4) S(X,X_3)}{n(n-1)^2} - \frac{g(X,X_3) S^2(X_2,X_4)}{(n-1)^2}.$$
Thus using \eqref{s2}, we get
$$[n S(X,X_1)-r g(X,X_1)] [n S(X_2,X_4)-r g(X_2,X_4)] = 0.$$
This implies that the manifold is Einstein. The converse part is obvious as on an Einstein manifold, $P=W$.
%-------------------------------
%
%===========================================================================
\begin{rem}
For a generalized curvature tensor $D$, $D \cdot \mathcal R = 0$ and $D\cdot R=0$ (resp., $D \cdot \mathcal S = 0$ and $D\cdot S=0$) give same structure but from Theorems \ref{thm4.5} (resp., Theorem \ref{thm4.6}) we can conclude that the structure $P \cdot \mathcal R = 0$ (resp., $P \cdot \mathcal S = 0$ and $P\cdot \mathcal P = 0$) is different from the structure $P\cdot R=0$ (resp., $P \cdot S = 0$ and $P\cdot P=0$). Similarly $P \cdot \mathcal R = L Q(g,\mathcal R)$ and $P\cdot R= L Q(g,R)$ (resp., $P \cdot S = L Q(g,S)$ and $P \cdot \mathcal S = L Q(g,\mathcal S)$, $P \cdot \mathcal P = L Q(g,\mathcal P)$ and $P\cdot P= L Q(g,P)$) give different structures.
\end{rem}
%=================================================================================

%%%%%%%%%%%%%%%%%%%%%%%%%%%%%%%%%%%%%%%%%%%%%%%%%%%%%%%%%%%%%%%%%%%%%%%%%%%%%%%%%%%%%%%%%%%%%%%%%%%%%%%%%%%%%%%%%%%%%%%%%
%                                P.R(0,4)=P.R(1,3), P.S(0,2)=P.S(1,1), P.P(0,4)=P.P(1,3)
%%%%%%%%%%%%%%%%%%%%%%%%%%%%%%%%%%%%%%%%%%%%%%%%%%%%%%%%%%%%%%%%%%%%%%%%%%%%%%%%%%%%%%%%%%%%%%%%%%%%%%%%%%%%%%%%%%%%%%%%%
\begin{pr}\label{pr4.4}\cite{DD91}
If $A$ be a symmetric (0,2)-tensor and $D$ be a generalized curvature tensor on a semi-Riemannian manifold $M$, then $Q(A,D) = 0$ implies either $A$ is of rank 1 or $B$ is linearly independent with $A\wedge A$.
\end{pr}
\indent From the above theorem we note that if $A$ is not of rank one then $Q(A,B)=0$ if and only if $B$ is linearly independent with $A\wedge A$. We also note that the result is not true for $D=P$. For the case of projective curvature tensor we have the following:
%=======================================================================
\begin{pr}\label{pr4.5}
For a symmetric (0,2) tensor $A$, if $Q(A,P)=0$ then $g$, $S$ and $A$ are linearly dependent.
\end{pr}
%-------------------------
\noindent \textbf{Proof:} From the condition $Q(A,P)=0$, we get
$$Q(A,R)=\frac{1}{n-1}\wedge_A\cdot\wedge_S.$$
Thus $\wedge_A\cdot\wedge_S$ possesses the following symmetry
$$\wedge_A\cdot\wedge_S(X_1,X_2,X_3,X_4,X,Y)=\wedge_A\cdot\wedge_S(X_3,X_4,X_1,X_2,X,Y).$$
Then taking contraction over $X_1$ and $X_4$ and using symmetry of $A$ and $S$, we get
\beb
&&\kappa A(X_3, Y) g(X, X_2) + \kappa A(X_2, Y) g(X, X_3) - \kappa A(X, X_3) g(X_2, Y) - \kappa A(X, X_2) g(X_3, Y)\\
&&- n A(X_3, Y) S(X, X_2) - n A(X_2, Y) S(X, X_3) + n A(X, X_3) S(X_2, Y) + n A(X, X_2) S(X_3, Y) = 0
\eeb
$$\Rightarrow Q(A,n S-\kappa g) =0.$$
Now by Lemma \ref{lem3.1}, $n S-\kappa g$ and $A$ are linearly dependent. Hence the result.
%-------------------------------
\begin{cor}\label{cor4.2}
Let $M$ be a semi-Riemannian manifold. If\\
(i) $Q(S,P)=0$ then $M$ is either Einstein or $\kappa = 0$.\\
(ii) $Q(g,P)=0$ then $M$ is Einstein.
\end{cor}
%
%%%%%%%%%%%%%%%%%%%%%%%%%%%%%%%%%%%%%%%%%%%%%%%%%%%%%%%%%%%%%%%%%%%%%%%%%%%%%%%%%%%%%%%%%%%%%%%%%%%%%%%%%%%%%%%%%%%%%
\begin{thm}\label{thm4.8}
If a semi-Riemannian manifold $M$ satisfies $P\cdot R = L Q(g,R)$, then\\
(i) $R\cdot R = L Q(g,R)$ if and only if $Q(S,R)=0$. Moreover if $S$ is not of rank 1, then $R\cdot R = L Q(g,R)$ if and only if $R=\lambda (S\wedge S)$ for some scalar $\lambda$.\\
(ii) $R\cdot R = 0$ if and only if $R=\lambda (\frac{1}{n-1}S+L g)\wedge(\frac{1}{n-1}S+L g)$ for some scalar $\lambda$ provided $M$ is not quasi-Einstein.\\
(iii) $P\cdot S = L Q(g,S)$ if and only if
$$R(Y,X_1,X_2,\mathcal S X) + R(Y,X_2,X_1,\mathcal S X) - R(X,X_1,X_2,\mathcal S Y) - R(X,X_2,X_1,\mathcal S Y) = 0.$$
(iv) $E = S^2 - \frac{n-1}{n-2} L(n S - \kappa g)$.
\end{thm}
%-------------------------
\noindent \textbf{Proof:}
(i) Since $P\cdot R = R\cdot R - \frac{1}{n-1} Q(S,R) = L Q(g,R)$, then $R\cdot R = L Q(g,R)$ if and only if $Q(S,R)=0$. Again if Rank($S$) is not equal to 1, then by using Proposition \ref{pr4.4}, 
\begin{center}
$(S,R)=0$ $\Leftrightarrow$ $R=\lambda (S\wedge S)$ for some scalar $\lambda$,\\
i.e., $R\cdot R = L Q(g,R)$ $\Leftrightarrow$ $R=\lambda (S\wedge S)$ for some scalar $\lambda$.\\
\end{center}
(ii) From the given hypothesis, 
\beb
&R\cdot R - \frac{1}{n-1}Q(S,R) = L Q(g, R)&\\
&\Rightarrow R\cdot R = Q(\frac{1}{n-1}S + L g, R).&
\eeb
Now if $M$ is not quasi-Einstein, then $(\frac{1}{n-1}S + L g)$ is not of rank one, Hence from Proposition \ref{pr4.4}, we get our assertion.\\
(iii) We know $R\cdot S = P\cdot S$, so contracting the condition $P\cdot R = R\cdot R - \frac{1}{n-1} Q(S,R) = L(g,R)$, we get 
\beb
&&R\cdot S[X_1, X_4, X,Y] - L Q(g,S) =\\\nonumber
&& -\frac{1}{n-1}\left[R(Y,X_1,X_4,\mathcal S X) + R(Y,X_4,X_1,\mathcal S X) - R(X,X_1,X_4,\mathcal S Y) - R(X,X_4,X_1,\mathcal S Y)\right].
\eeb
This implies (iii).\\
(iv) Again contracting the above equation over $X_1$ and $X_4$, we get (iv).\\
%=================================
\indent We know from the definitions that $P\cdot R - L Q(S,R) = R\cdot R - (L+\frac{1}{n-1}) Q(S,R)$. Hence using the same process of the proof of the Theorem \ref{thm4.9}, we get the following:
\begin{thm}\label{thm4.9}
If a semi-Riemannian manifold $M$ satisfies $P\cdot R = L Q(S,R)$, then\\
(i) $R\cdot R = (L+\frac{1}{n-1}) Q(S,R)$.\\
(ii) $R\cdot R = 0$ if and only if $Q(S,R)=0$ or $L=-\frac{1}{n-1}$. Moreover if $M$ is Ricci simple, then $P\cdot R = 0$ if and only if $R=\lambda (S\wedge S)$ for some scalar $\lambda$.\\
(iii) $P\cdot S = 0$ if and only if
$$\left(L+\frac{1}{n-1}\right)\left[R(Y,X_1,X_2,\mathcal S X) + R(Y,X_2,X_1,\mathcal S X) - R(X,X_1,X_2,\mathcal S Y) - R(X,X_2,X_1,\mathcal S Y)\right] = 0.$$
(iv) $\left(L-\frac{n-2}{n-1}\right)(E - S^2) = 0$.
\end{thm}
%=================================
\begin{cor}
Let $M$ be a semi-Riemannian manifold $M$ satisfying $P\cdot R =0$. Then\\
(i) $R\cdot R =0$ if and only if $Q(S,R)=0$. Moreover if $S$ is not Ricci simple, then $R\cdot R =0$ if and only if 
$R=\lambda (S\wedge S)$ for some scalar $\lambda$.\\
(ii) $P\cdot S =0$ if and only if
\be\label{cond}
R(Y,X_1,X_2,\mathcal S X) + R(Y,X_2,X_1,\mathcal S X) - R(X,X_1,X_2,\mathcal S Y) - R(X,X_2,X_1,\mathcal S Y) = 0.
\ee
(iii) $E = S^2$.
\end{cor}
%%%%%%%%%%%%%%%%%%%%%%%%%%%%%%%%%%%%%%%%%%%%%%%%%%%%%%%%%%%%%%%%%%%%%%%%%%%%%%%%%%%%%%%%%%%%%%%%%%%%%%%%%%%%%%%%%%%%%
%==============================================================================
\begin{thm}\label{thm4.10}
If a semi-Riemannian manifold $M$ satisfies 
\be\label{p.pg}
(P\cdot P)(X_1,X_2,X_3,X_4,X,Y) = L Q(g,P)(X_1,X_2,X_3,X_4,X,Y),
\ee
then\\
(i) $P\cdot R = L Q(g,R)$ if and only if $R\cdot (\wedge_S) = L Q(g,\wedge_S) + \frac{1}{n-1} Q(S,\wedge_S)$.\\
(ii) $R\cdot P = L Q(g,P)$ if and only if $Q(S,P) = 0$.\\
(iii) $R\cdot R = L Q(g,R)$ if and only if $Q(S,P) = 0$.\\
(iv) $n(n-1) R\cdot S = (n^2 L - n L + \kappa) Q(g,S)$.\\
(v) $(n-1)E = \kappa S - S^2$,\\
(vi) $n^2 S^2 - 2 n \kappa S + \kappa^2 g = 0$ and $\kappa^{(2)} = \frac{r^2}{n}$. Moreover if $M$ is a Riemannian manifold, then it is an Einstein manifold.\\
(vii) $L n^2 (n S - \kappa g) = 0$. Moreover if $L$ is nowhere zero, then $M$ is an Einstein manifold.
\end{thm}
%-----------------------
\noindent \textbf{Proof:} (i) Since $P = R-\frac{1}{n-1}\wedge_S$, so 
\beb
[P\cdot P - L Q(g,P)] - [P\cdot R - L Q(g,R)] &=& [P\cdot R - \frac{1}{n-1}P\cdot(\wedge_S) - L Q(g,R) + \frac{L}{n-1}Q(g,\wedge_S)]\\
&&- [P\cdot R - L Q(g,R)]\\
&=& - \frac{1}{n-1}[P\cdot(\wedge_S) - \frac{L}{n-1}Q(g,\wedge_S)]
\eeb
Thus $M$ satisfies $P\cdot R = L Q(g,R)$ if and only if $P\cdot(\wedge_S) = \frac{L}{n-1}Q(g,\wedge_S)$, i.e., 
$R\cdot(\wedge_S) = \frac{1}{n-1}Q(S,\wedge_S) + \frac{L}{n-1}Q(g,\wedge_S)$.\\
(ii) Again
\beb
[P\cdot P - L Q(g,P)] - [R\cdot P - L Q(g,P)] &=& [R\cdot P - \frac{1}{n-1}(\wedge_S)\cdot P - L Q(g,P)]\\
&&- [R\cdot P - L Q(g,P)]\\
&=& - \frac{1}{n-1}(\wedge_S)\cdot P
\eeb
Thus $M$ satisfies $R\cdot P = L Q(g,P)$ if and only if $(\wedge_S)\cdot P = 0$, i.e., $Q(S,P) = 0$.\\
(iii) From Corollary 6.3 of \cite{SK14}, we know that the curvature conditions $R\cdot P = L Q(g,P)$ and $R\cdot R = L Q(g,R)$ are equivalent. Hence $M$ satisfies $R\cdot R = L Q(g,R)$ if and only if $Q(S,P) = 0$.\\
(iv) Since $P(X_1,X_2,X_3,X_4) \neq P(X_3,X_4,X_1,X_2)$, then contracting 
\beb
&&[(P\cdot P)(X_1,X_2,X_3,X_4,X,Y) - L Q(g,P)(X_1,X_2,X_3,X_4,X,Y)]\\
&&-[(P\cdot P)(X_3,X_4,X_1,X_2,X,Y) - L Q(g,P)(X_3,X_4,X_1,X_2,X,Y)] = 0,
\eeb
over $X_1$ and $X_4$, we get
\beb
&&-\frac{L n S(X,X_2) g(Y,X_3)}{n-1}-\frac{L n S(X,X_3) g(Y,X_2)}{n-1}+\frac{L n g(X,X_3) S(Y,X_2)}{n-1}\\
&&+\frac{L n g(X,X_2) S(Y,X_3)}{n-1}-\frac{\kappa S(X,X_2) g(Y,X_3)}{(n-1)^2}-\frac{\kappa S(X,X_3) g(Y,X_2)}{(n-1)^2}+\frac{\kappa g(X,X_3) S(Y,X_2)}{(n-1)^2}\\
&&+\frac{\kappa g(X,X_2) S(Y,X_3)}{(n-1)^2}+\frac{n R(X,Y,X_2,S(X_3))}{n-1}+\frac{n R(X,Y,X_3,S(X_2))}{n-1} = 0.
\eeb
Now rearranging the above equation, we get
\beb
&\left(\frac{L n}{n-1}+\frac{r}{(n-1)^2}\right) Q(g,S)(X_2,X_3,X,Y)-\frac{n}{n-1} (R\cdot S)(X_2,X_3,X,Y) = 0&\\
&\Rightarrow n(n-1) R\cdot S = (n^2 L - n L + \kappa) Q(g,S).&
\eeb
(v) Now contracting the condition \eqref{p.pg} over $X_2$ and $X_4$, we get
\beb
&&\frac{S^2(X,X_3) g(Y,X_1)}{(n-1)^2}-\frac{g(X,X_1) S^2(Y,X_3)}{(n-1)^2}-\frac{n R(X,X_1,X_3,S(Y))}{(n-1)^2}\\
&&+\frac{R(X,X_1,X_3,S(Y))}{(n-1)^2}-\frac{n R(X,X_3,X_1,S(Y))}{(n-1)^2}+\frac{R(X,X_3,X_1,S(Y))}{(n-1)^2}\\
&&-\frac{n R(X_1,Y,X_3,S(X))}{(n-1)^2}+\frac{R(X_1,Y,X_3,S(X))}{(n-1)^2}-\frac{n R(X_1,S(X),X_3,Y)}{(n-1)^2}\\
&&+\frac{R(X_1,S(X),X_3,Y)}{(n-1)^2}-\frac{S(X,X_3)S(Y,X_1)}{(n-1)^2}+\frac{S(X,X_1) S(Y,X_3)}{(n-1)^2} = 0.
\eeb
Again contracting the above over $Y$ and $X_1$, we get 
$$\frac{1}{(n-1)^2}\left[(n-1) E(X,X_3) - \kappa S(X,X_3) + S^2(X,X_3)\right] = 0 \ \Rightarrow \ (n-1)E = \kappa S - S^2.$$
(vi) Since $P(X_1,X_2,X_3,X_4) \neq - P(X_1,X_2,X_4,X_3)$, contracting the condition \eqref{p.pg} over $X_3$ and $X_4$, we get
\beb
&&-\frac{S^2(X,X_1) g(Y,X_2)}{(n-1)^2}+\frac{S^2(X,X_2) g(Y,X_1)}{(n-1)^2}+\frac{g(X,X_2) S^2(Y,X_1)}{(n-1)^2}\\
&&-\frac{g(X,X_1) S^2(Y,X_2)}{(n-1)^2}-\frac{2 S(X,X_2) S(Y,X_1)}{(n-1)^2}+\frac{2 S(X,X_1)S(Y,X_2)}{(n-1)^2} = 0.
\eeb
Again contracting the above over $Y$ and $X_1$, we get
$$\frac{1}{(n-1)^2}[\kappa^{(2)} g(X,X_2)+n S^2(X,X_2)-2 \kappa S(X,X_2)] = 0,$$
which implies $\kappa^{(2)} = \frac{\kappa^2}{n}$, and hence from above equation we say that $n^2 S^2 - 2 n \kappa S + \kappa^2 g = 0$. The next part directly follows from Proposition \ref{pr4.6}.\\
(vii) Now contracting \eqref{p.pg} over $X_2$ and $X_3$ and then replacing $X_4$ by $X_3$, we get
\beb
&&-\frac{n R(X,Y,X_1,S(X_3))}{n-1}-\frac{n R(X,Y,X_3,S(X_1))}{n-1}+\frac{R(X,X_1,X_3,S(Y))}{n-1}+\frac{R(X,X_3,X_1,S(Y))}{n-1}\\
&&+\frac{R(X_1,Y,X_3,S(X))}{n-1}+\frac{R(X_1,S(X),X_3,Y)}{n-1}-\frac{(L (n-1) n+\kappa)}{(n-1)^2} Q(g, S)(X_1, X_3, X, Y)\\
&&-\frac{S^2(X,X_1) g(Y,X_3)}{(n-1)^2}+\frac{g(X,X_3) S^2(Y,X_1)}{(n-1)^2}-\frac{S(X,X_3) S(Y,X_1)}{(n-1)^2}+\frac{S(X,X_1)
   S(Y,X_3)}{(n-1)^2} = 0
\eeb
Again contracting the above over $Y$ and $X_1$, we get
$$E(X,X_3) = -\frac{\kappa^{(2)}-\kappa (L (n-1) n+\kappa)}{(n-1)^2} g(X,X_3) - \frac{L n^2+\kappa}{n-1} S(X,X_3) + S^2(X,X_3)$$
$$\Rightarrow E = -\frac{\kappa^{(2)}-\kappa (L (n-1) n+\kappa)}{(n-1)^2} g - \frac{L n^2+\kappa}{n-1} S + S^2.$$
Now using (v), we get
$$\frac{\kappa}{n-1} S - \frac{1}{n-1} S^2 = -\frac{\kappa^{(2)}-\kappa (L (n-1) n+\kappa)}{(n-1)^2} g - \frac{L n^2+\kappa}{n-1} S + S^2$$
$$\Rightarrow (-L (-1 + n) n \kappa - \kappa^2 + \kappa^{(2)}) g +  (n-1)(L n^2 + 2 \kappa) S - n(n-1) S^2 = 0.$$
Again using (vi), we get $L n (n S - \kappa g) = 0$. This completes the proof.\\
%%%%%%%%%%%%%%%%%%%%%%%%%%%%%%%%%%%%%%%%%%%%%%%%%%%%%%%%%%%%%%%%%%%%%%%%%%%%%%%%%%%%%%%
\indent Using the same technique of the proof of the previous theorem, we get the following:
\begin{thm}\label{thm4.11}
If a semi-Riemannian manifold $M$ satisfies 
\be\label{p.ps}
(P\cdot P)(X_1,X_2,X_3,X_4,X,Y) = L Q(S,P)(X_1,X_2,X_3,X_4,X,Y),
\ee
then\\
(i) $P\cdot R = L Q(S,R)$ if and only if $R\cdot (\wedge_S) = \left(L + \frac{1}{n-1}\right) Q(S,\wedge_S)$.\\
(ii) $R\cdot R = L Q(S,R)$ if and only if $R\cdot (\wedge_S) + Q(S,R) = \left(L + \frac{1}{n-1}\right) Q(S,\wedge_S)$.\\
(iii) $R\cdot P = L Q(S,P)$ if and only if $Q(S,P) = 0$.\\
(iv) $n(n-1) R\cdot S = (1+(n-1)\kappa L) Q(g,S)$.\\
(v) $(n-1)E = \kappa S - S^2$,\\
(vi) $n^2 S^2 - 2 n \kappa S + \kappa^2 g = 0$ and $\kappa^{(2)} = \frac{r^2}{n}$ if $L \neq -\frac{1}{n-1}$. Moreover if $M$ is a Riemannian manifold and $L \neq -\frac{1}{n-1}$, then it is an Einstein manifold.\\
(vii) $L \kappa (n S - \kappa g) = 0$.
\end{thm}
%%%%%%%%%%%%%%%
\begin{cor}\label{cor4.4}
If a semi-Riemannian manifold $M$ satisfies $P\cdot P =0$, then\\
(i) $P\cdot R = 0$ if $P\cdot(\wedge_S) = 0$.\\
(ii) $R\cdot P = 0$ or $R\cdot R = 0$ if $Q(S,P) = 0$.\\
(iii) $(n-1)E = \kappa S - S^2$,\\
(iv) $n^2 S^2 - 2 n \kappa S + \kappa^2 g = 0$ and thus $\kappa^{(2)} = \frac{r^2}{n}$, where $\kappa^{(2)}$ is the trace of $S^2$. Moreover if $M$ is a Riemannian manifold, then it is an Einstein manifold.\\
(v) $n(n-1)R\cdot S = \kappa Q(g, S)$.\\
(vi) $S\wedge S = g\wedge S^2$.
\end{cor}
%=========================================
%=====================================================================================================
%=== === === === === === === === === === === ===  P.R (1, 3) and P.S (1, 1) case === === === === === == = === === === 
%===================================================================================================
\indent Now we can easily check that the tensor $g((P(X,Y)\cdot \mathcal S)(X_1), X_2)$ is not symmetric in $X_1$ and $X_2$. Now using this asymmetry, we get the following results:
\begin{thm}\label{thm4.12}
If a semi-Riemannian manifold $M$ satisfies $P\cdot \mathcal S = L Q(g,\mathcal S)$, then\\
(i) $S \wedge S = g\wedge S^2$.\\
(ii) $n^2 S^2 - 2 n \kappa S + \kappa^2 g = 0$ and $\kappa^{(2)} = \frac{r^2}{n}$. Moreover if $M$ is a Riemannian manifold, then $M$ is an Einstein manifold.\\
(iii) $(n-1)E = L (n-1) \kappa g +(\kappa-L (n-1) n) S - S^2$.
\end{thm}
%-------------------------
\noindent \textbf{Proof:} (i) Since $P\cdot \mathcal S = L Q(g,\mathcal S)$, so
\beb
g(P\cdot \mathcal S(X_1,X,Y),X_2) &-& L g(Q(g,\mathcal S)(X_1,X,Y),X_2)\\
 &=& g(P\cdot \mathcal S(X_2,X,Y),X_1) - L g(Q(g,\mathcal S)(X_2,X,Y),X_1)
\eeb
\beb
\Rightarrow && \frac{S^2(X,X_1) g(X_2,Y)}{n-1}-\frac{S^2(X,X_2) g(X_1,Y)}{n-1}-\frac{g(X,X_2) S^2(X_1,Y)}{n-1}\\
						&&+\frac{g(X,X_1)S^2(X_2,Y)}{n-1}+\frac{2 S(X,X_2) S(X_1,Y)}{n-1}-\frac{2 S(X,X_1)S(X_2,Y)}{n-1} = 0.
\eeb
$$\Rightarrow \frac{1}{n-1}\left[(S\wedge S)(X_1,X_2,X,Y)-(g\wedge S^2)(X_1,X_2,X,Y)\right] = 0$$
(ii) Now contracting the above equation over $X$ and $X_1$, we get
$$\frac{\kappa^{(2)} g(X_2,Y)}{n-1}-\frac{2 \kappa S(X_2,Y)}{n-1}+\frac{n S^2(X_2,Y)}{n-1} = 0,$$
which implies $\kappa^{(2)} = \frac{1}{n}\kappa$, and hence $n^2 S^2 - 2 n \kappa S + \kappa^2 g = 0$.\\
(iii) Now contracting the given condition $g((P\cdot \mathcal S)(X_1,X,Y),X_2) = L g(Q(g,\mathcal S)(X_1,X,Y),X_2)$ over $X$ and $X_2$, we get
$$(n-1)E(Y,X_2) = L (n-1) \kappa g(Y,X_1)+(\kappa-L (n-1) n) S(Y,X_1)-S^2(Y,X_1).$$
%%============================================================================================================
\begin{cor}
If a semi-Riemannian manifold $M$ satisfies $P\cdot \mathcal S =0$, then\\
(i) $S \wedge S = g\wedge S^2$.\\
(ii) $n^2 S^2 - 2 n \kappa S + \kappa^2 g = 0$ and $\kappa^{(2)} = \frac{r^2}{n}$. Moreover if $M$ is a Riemannian manifold, then $M$ is an Einstein manifold.\\
(iii) $(n-1)E = \kappa S - S^2.$
\end{cor}
%=============================================================================================================
%=============================================================================================================
\begin{thm}\label{thm4.13}
If a semi-Riemannian manifold $M$ satisfies $P\cdot \mathcal R = L Q(g, \mathcal R)$, then\\
(i) $P\cdot S = R\cdot S = L Q(g,S)$, i.e., $M$ is Ricci pseudosymmetric and thus $E = L \kappa g - L n S + S^2$.\\
(ii) $n^2 S^2 - 2 n \kappa S + \kappa^2 g = 0$ and $\kappa^{(2)} = \frac{r^2}{n}$. Moreover if $M$ is a Riemannian manifold, then $M$ is an Einstein manifold.\\
(iii) $S\wedge S = g\wedge S^2$.\\
(iv) $P\cdot \mathcal S = L Q(g, \mathcal S)$ if and only if \eqref{cond} holds or
$$S(X, X_2) S(X_1, Y) - S(X, X_1) S(X_2, Y) +  g(X_2, Y) S^2(X, X_1) - g(X, X_2) S^2(X_1, Y) =0.$$
Especially, if $Q(S,R)=0$ on $M$, then $P\cdot \mathcal R =0$ $\Rightarrow$ $P\cdot \mathcal S =0$.
\end{thm}
%-------------------------
\noindent \textbf{Proof:} From the given hypothesis
\beb
&&\mathcal P(X,Y)\mathcal R(X_1, X_2)X_3 - \mathcal R(\mathcal P(X,Y)X_1, X_2)X_3\\\nonumber
&&-\mathcal R(X_1, \mathcal P(X,Y)X_2)X_3 - \mathcal R(X_1, X_2)\mathcal P(X,Y)X_3\\\nonumber
&&\hspace{1in} = L(X \wedge Y)\mathcal R(X_1, X_2)X_3 - L \mathcal R((X \wedge Y)X_1, X_2)X_3\\\nonumber
&&\hspace{1in} -L \mathcal R(X_1, (X \wedge Y)X_2)X_3 - L \mathcal R(X_1, X_2)(X \wedge Y)X_3.
\eeb
 Since $g$ is non-degenerate, the above condition is equivalent to
\bea\label{p.rg13}
&&P(X,Y,\mathcal R(X_1, X_2)X_3,X_4) - R(\mathcal P(X,Y)X_1, X_2,X_3,X_4)\\\nonumber
&&-R(X_1, \mathcal P(X,Y)X_2,X_3,X_4) - R(X_1, X_2,\mathcal P(X,Y)X_3,X_4)\\\nonumber
&&\hspace{1cm} = L (X\wedge Y)(\mathcal R(X_1, X_2)X_3,X_4) - L R((X \wedge Y)X_1, X_2,X_3,X_4)\\\nonumber
&&\hspace{1cm} -L R(X_1, (X \wedge Y)X_2,X_3,X_4) -L R(X_1, X_2,(X \wedge Y)X_3,X_4).
\eea
(i) Taking contraction over $X_1$ and $X_4$ in \eqref{p.rg13}, we get 
\beb
&&-R(X,Y,X_2,S(X_3))-R(X,Y,X_3,S(X_2))\\
&&+L (S(X,X_2) g(Y,X_3)+S(X,X_3) g(Y,X_2)-g(X,X_3) S(Y,X_2)-g(X,X_2)S(Y,X_3)) = 0
\eeb
i.e., $P\cdot S = R\cdot S = LQ(g,S)$. Again contacting this we get $E = L \kappa g - L n S + S^2$.\\
(ii) Contracting \eqref{p.rg13} over $X_2$ and $X_3$, we get 
\bea\label{p.rg1323}
&&L S(X,X_1) g(Y,X_4)+L S(X,X_4) g(Y,X_1)-L g(X,X_4) S(Y,X_1)-L g(X,X_1) S(Y,X_4)\\\nonumber
&&+\frac{S^2(X,X_1)g(Y,X_4)}{n-1}-\frac{g(X,X_4)S^2(Y,X_1)}{n-1}+\frac{R(X,X_1,X_4,S(Y))}{n-1}\\\nonumber
&&+\frac{R(X,X_4,X_1,S(Y))}{n-1}+\frac{R(X_1,Y,X_4,S(X))}{n-1}+\frac{R(X_1,S(X),X_4,Y)}{n-1}\\\nonumber
&&+\frac{S(X,X_4)S(Y,X_1)}{n-1}-\frac{S(X,X_1)S(Y,X_4)}{n-1}-R(X,Y,X_1,S(X_4))-R(X,Y,X_4,S(X_1)) =0.
\eea
Now contracting \eqref{p.rg1323} over $Y$ and $X_4$, we get
$$\frac{-L \kappa g(X,X_1)+S(X,X_1) (L n-\kappa)+n S^2(X,X_1)}{n-1} = 0.$$
Again contracting \eqref{p.rg1323} over $Y$ and $X_4$ and then replacing $X_4$ by $X_1$, we get
$$\frac{g(X,X_1) (-L \kappa-\kappa^{(2)})+S(X,X_1) (L n+\kappa)}{n-1} = 0.$$
Now from last two equation we get $\kappa^{(2)} g + n S^2 -2 \kappa S = 0$, which implies $n^2 S^2 - 2 n \kappa S + \kappa^2 g = 0$ and $\kappa^{(2)} = \frac{r^2}{n}$.\\
(iii) Interchanging $X_4$ and $X_1$ in \eqref{p.rg1323}, we get
\beb
&&L S(X,X_4) g(Y,X_1)+L S(X,X_1) g(Y,X_4)-L g(X,X_1) S(Y,X_4)-L g(X,X_4) S(Y,X_1)\\\nonumber
&&+\frac{S^2(X,X_4)g(Y,X_1)}{n-1}-\frac{g(X,X_1)S^2(Y,X_4)}{n-1}+\frac{R(X,X_4,X_1,S(Y))}{n-1}\\\nonumber
&&+\frac{R(X,X_1,X_4,S(Y))}{n-1}+\frac{R(X_4,Y,X_1,S(X))}{n-1}+\frac{R(X_4,S(X),X_1,Y)}{n-1}\\\nonumber
&&+\frac{S(X,X_1)S(Y,X_4)}{n-1}-\frac{S(X,X_4)S(Y,X_1)}{n-1}-R(X,Y,X_4,S(X_1))-R(X,Y,X_1,S(X_4)) =0.
\eeb
Now subtracting the above equation from \eqref{p.rg1323}, we get 
\beb
&&\frac{S^2(X,X_1) g(Y,X_4)}{n-1}-\frac{S^2(X,X_4) g(Y,X_1)}{n-1}-\frac{g(X,X_4) S^2(Y,X_1)}{n-1}\\
&&+\frac{g(X,X_1)S^2(Y,X_4)}{n-1}+\frac{2 S(X,X_4) S(Y,X_1)}{n-1}-\frac{2 S(X,X_1)S(Y,X_4)}{n-1} = 0.
\eeb
$$\Rightarrow \frac{1}{n-1}[(S\wedge S)(X_1,X_4,X,Y) - (g\wedge S^2)(X_1,X_4,X,Y)] = 0.$$
This completes the proof.\\
(iv) Now $P\cdot \mathcal S = L Q(g, \mathcal S)$ holds if and only if
\beb
&&L S(X,X_1) g(X_4,Y)+L S(X,X_4) g(X_1,Y)-L g(X,X_4) S(X_1,Y)-L g(X,X_1) S(X_4,Y)\\
&&+\frac{S^2(X,X_1)g(X_4,Y)}{n-1}-\frac{g(X,X_4) S^2(X_1,Y)}{n-1}+\frac{S(X,X_4)S(X_1,Y)}{n-1}\\
&&-\frac{S(X,X_1)S(X_4,Y)}{n-1}-R(X,Y,X_1,S(X_4))-R(X,Y,X_4,S(X_1)) = 0.
\eeb
Hence from \eqref{p.rg1323}, we say that $P\cdot \mathcal S = L Q(g, \mathcal S)$ holds if and only if
\beb
\frac{1}{n-1}[R(X,X_1,X_4,S(Y))+R(X,X_4,X_1,S(Y))+R(X_1,Y,X_4,S(X))+R(X_1,S(X),X_4,Y)]=0.
\eeb
%===========================================================================================================
\begin{thm}\label{thm4.14}
If a semi-Riemannian manifold $M$ satisfies $P\cdot \mathcal R = L Q(S, \mathcal R)$, then\\
(i) $R\cdot R = 0$ if $L+\frac{1}{n-1} = 0$.\\
(ii) $P\cdot S = R\cdot S = 0$, i.e., $M$ is Ricci semisymmetric and thus $E = S^2$.\\
(iii) $n^2 S^2 - 2 n \kappa S + \kappa^2 g = 0$ and $\kappa^{(2)} = \frac{r^2}{n}$ if $L+\frac{1}{n-1}\neq 0$.\\
(iv) $S\wedge S = g\wedge S^2$.\\
$S^2 = \frac{\kappa}{n} S = \frac{\kappa^{(2)}}{n} g$. Moreover if $\kappa \ne 0$, then $M$ is Einstein.\\
(v) for $L+\frac{1}{n-1}\neq 0$, $P\cdot \mathcal S = L Q(g, \mathcal S)$ if and only if \eqref{cond} holds.
\end{thm}
%-------------------------
\noindent \textbf{Proof:} From the given hypothesis
\beb
&&\mathcal P(X,Y)\mathcal R(X_1, X_2)X_3 - \mathcal R(\mathcal P(X,Y)X_1, X_2)X_3\\\nonumber
&&-\mathcal R(X_1, \mathcal P(X,Y)X_2)X_3 - \mathcal R(X_1, X_2)\mathcal P(X,Y)X_3\\\nonumber
&&\hspace{1in} =L (X \wedge_S Y)\mathcal R(X_1, X_2)X_3 - L \mathcal R((X \wedge_S Y)X_1, X_2)X_3\\\nonumber
&&\hspace{1in} -L \mathcal R(X_1, (X \wedge_S Y)X_2)X_3 - L \mathcal R(X_1, X_2)(X \wedge_S Y)X_3.
\eeb
 Since $g$ is non-degenerate, the above condition is equivalent to
\bea\label{p.rs}
&&P(X,Y,\mathcal R(X_1, X_2)X_3,X_4) - R(\mathcal P(X,Y)X_1, X_2,X_3,X_4)\\\nonumber
&&-R(X_1, \mathcal P(X,Y)X_2,X_3,X_4) - R(X_1, X_2,\mathcal P(X,Y)X_3,X_4)\\\nonumber
&&\hspace{1cm} = L (X\wedge_S Y)(\mathcal R(X_1, X_2)X_3,X_4) - L R((X \wedge_S Y)X_1, X_2,X_3,X_4)\\\nonumber
&&\hspace{1cm} - L R(X_1, (X \wedge_S Y)X_2,X_3,X_4) - L R(X_1, X_2,(X \wedge_S Y)X_3,X_4).
\eea
(i) The proof is obvious, since $P=R-\frac{1}{n-1}\wedge_S$.\\
(ii) Contracting \eqref{p.rs} over $X_1$ and $X_4$ we get the result.\\
(iii) Contracting \eqref{p.rs} over $X_2$ and $X_3$, we get
\bea\label{p.rs23}
&&(L+\frac{1}{n-1}) [S^2(X,X_1) g(Y,X_4)-g(X,X_4) S^2(Y,X_1)\\\nonumber
											&&\hspace{3cm} +S(X,X_4) S(Y,X_1)-S(X,X_1)S(Y,X_4)]\\\nonumber
&&+(L+\frac{1}{n-1}) [R(X_1,Y,X_4,S(X))+R(X_1,S(X),X_4,Y)\\\nonumber
											&&\hspace{3cm} +R(X,X_1,X_4,S(Y))+R(X,X_4,X_1,S(Y))]\\\nonumber
&&-R(X,Y,X_1,S(X_4))-R(X,Y,X_4,S(X_1)) = 0
\eea
Now contracting \eqref{p.rs23} over $X_4$ and $Y$ and putting $E = S^2$, we get
$$(L+\frac{1}{n-1})(n S^2 - \kappa S) = 0.$$
Again contracting \eqref{p.rs23} over $X_1$ and $Y$ and putting $E = S^2$, we get
$$(L+\frac{1}{n-1})(\kappa S - \kappa^{(2)} g) = 0.$$
Thus from the last two equations we get $(L+\frac{1}{n-1})(\kappa^{(2)} g - 2 \kappa S + n S^2) =0$. Hence $n^2 S^2 - 2 n \kappa S + \kappa^2 g = 0$ and $\kappa^{(2)} = \frac{r^2}{n}$ if $L+\frac{1}{n-1}\neq 0$\\
The proof of (iv) and (v) are similar to the proof of (iii) and (iv) of Theorem \ref{thm4.13}. 
%================================================================================
\begin{cor}
If a semi-Riemannian manifold $M$ satisfies $P\cdot \mathcal R =0$, then\\
(i) $P\cdot S = R\cdot S = 0$, i.e., $M$ is Ricci semisymmetric and thus $E = S^2$\\
(ii) $S^2 = \frac{\kappa}{n} S = \frac{\kappa^{(2)}}{n} g$. Moreover if $\kappa \ne 0$, then $M$ is Einstein.\\
(iii) $n^2 S^2 - 2 n \kappa S + \kappa^2 g = 0$ and $\kappa^{(2)} = \frac{r^2}{n}$. Moreover if $M$ is a Riemannian manifold, then $M$ is an Einstein manifold.\\
(iv) $P\cdot \mathcal S =0$ if and only if \eqref{cond} holds or
$$S(X, X_2) S(X_1, Y) - S(X, X_1) S(X_2, Y) +  g(X_2, Y) S^2(X, X_1) - g(X, X_2) S^2(X_1, Y) =0.$$
Especially, if $Q(S,R)=0$ on $M$, then $P\cdot \mathcal R =0$ $\Rightarrow$ $P\cdot \mathcal S =0$.
\end{cor}
%-------------------------
%=====================================================================
Now by similar technique of the proof of the Theorem \ref{thm4.13} and \ref{thm4.14}, we get the following theorem.
%=====================================================================
\begin{thm}\label{thm4.15}
If a semi-Riemannian manifold $M$ satisfies $P\cdot \mathcal P = L Q(g, \mathcal P)$, then\\
(i) $(n-1)E = L n \kappa g - (L n^2 + \kappa) S - S^2$,\\
(ii) $n^2 S^2 - 2 n \kappa S + \kappa^2 g = 0$ and thus $\kappa^{(2)} = \frac{r^2}{n}$.\\
(iii) $\frac{n}{n - 1} (R\cdot \mathcal S) = \frac{1}{n - 1} Q(S,\mathcal S) - \frac{n L}{n - 1} Q(g,\mathcal S)$ holds if \eqref{cond} holds.\\
(iv) $S\wedge S = g\wedge S^2$.
\end{thm}
%============================================
\begin{thm}\label{thm4.16}
If a semi-Riemannian manifold $M$ satisfies $P\cdot \mathcal P = L Q(S, \mathcal P)$, then\\
(i) $(L-1)(n-1)E = -(1 + L (-1 + n)) \kappa S + (1 + L (-1 + n^2)) S^2$,\\
(ii) $n^2 S^2 - 2 n \kappa S + \kappa^2 g = 0$ and $\kappa^{(2)} = \frac{r^2}{n}$ if $L + \frac{1}{n-1} \neq 0$.\\
(iii) $\frac{n}{n - 1} (P\cdot \mathcal S) = \left(L-\frac{1}{(n - 1)^2}\right) Q(S,\mathcal S)$ holds if \eqref{cond} holds and $L + \frac{1}{n-1} \neq 0$.\\
(iv) $S\wedge S = g\wedge S^2$.
\end{thm}
%============================================
\begin{cor}
If a semi-Riemannian manifold $M$ satisfies $P\cdot \mathcal P =0$, then\\
(i) $(n-1)E = \kappa S - S^2$,\\
(ii) $n^2 S^2 - 2 n \kappa S + \kappa^2 g = 0$ and hence $\kappa^{(2)} = \frac{\kappa^2}{n}$.\\
(iii) $\frac{n}{n - 1} (R\cdot \mathcal S) = \frac{1}{n - 1} Q(S,\mathcal S)$ holds if \eqref{cond} holds.\\
(iv) $S\wedge S = g\wedge S^2$.
\end{cor}
%%%%%%%%%%%%%%%%%%%%%%%%%%%%%%%%%%%%%%%%%%%%%%%%%%%%%%%%%%%%%%%%%%%%%%%%%%%%%%%%%%%%%%%%%%%%%%%%%%%%%%%%%%%%%%%%%%%%%%%%%
%                                                  Roter type cases
%%%%%%%%%%%%%%%%%%%%%%%%%%%%%%%%%%%%%%%%%%%%%%%%%%%%%%%%%%%%%%%%%%%%%%%%%%%%%%%%%%%%%%%%%%%%%%%%%%%%%%%%%%%%%%%%%%%%%%%%%
\begin{lem}\label{lem5.1}\cite{SKgrt}
A Roter type semi-Riemannian manifold satisfies 
$$2\alpha S^2 = (\kappa \beta + 2 (-1 + n) \gamma) g + (2 \kappa\alpha + (-2 + n)\beta) S.$$
\end{lem}
\begin{lem}\label{lem5.2}\cite{SKgrt}
A generalized Roter type semi-Riemannian manifold is of constant curvature if and only if it is Einstein.
\end{lem}
%============================
Now using Theorem \ref{thm4.8}, Theorem \ref{thm4.9} and Lemma \ref{lem5.2}, we get the following generalization of the main result (Theorem 1) of \cite{DRV89} for non-conformally flat case. 
\begin{thm}
On a generalized Roter type Riemannian manifold, the following conditions are equivalent:\\
(i) $P\cdot \mathcal S = 0,$ \ \ (ii) $P\cdot \mathcal S = L Q(g,\mathcal S),$\\
(iii) $P\cdot \mathcal R = 0,$ \ \ (iv) $P\cdot \mathcal R = L Q(g,\mathcal R)$,  \ \ (v) $P\cdot \mathcal R = L Q(S,\mathcal R),$ $L\neq \frac{-1}{n-1}$,\\
(vi) $P\cdot P = 0,$ \ \ (vii) $P\cdot P = L Q(g,P)$,  \ \ (viii) $P\cdot P = L Q(S,P),$ $L\neq \frac{-1}{n-1}$,\\
(ix) $P\cdot \mathcal P = 0,$ \ \ (x) $P\cdot \mathcal P = L Q(g,\mathcal P)$,  \ \ (xi) $P\cdot \mathcal P = L Q(S,\mathcal P),$ $L\neq \frac{-1}{n-1}$,\\
(xii) $M$ is a manifold of constant curvature, where $L$ is a scalar.
\end{thm}
%%%%%%%%%%%%%%%%%%%%%%%%%%%%%%%%%%%%%%%%%%%%%%%%%%%%%%%%%%%%%%%%%%%%%%%%%%%%%%%%%%%%%%%%%%%%%%%%%%%%%%%%%%%%%%%%%%%%%%%%%
%%%%%%%%%%%%%%%%%%%%%%%%%%%%%%%%%%%%%%%%%%%%%%%%%%%%%%%%%%%%%%%%%%%%%%%%%%%%%%%%%%%%%%%%%%%%%%%%%%%%%%%%%%%%%%%%%%%%%%%%%
%                                              Examples
%%%%%%%%%%%%%%%%%%%%%%%%%%%%%%%%%%%%%%%%%%%%%%%%%%%%%%%%%%%%%%%%%%%%%%%%%%%%%%%%%%%%%%%%%%%%%%%%%%%%%%%%%%%%%%%%%%%%%%%%%
\section{\bf Examples}\label{exam}
%%%%%%%%%%%%%%%%%%%%%%%%%%%%%%%%%%
\textbf{Example 1:} Let $M_1$ be a $4$-dimensional connected semi-Riemannian manifold endowed with the semi-Riemannian metric
\be\label{met1}
ds^2 = e^{x^1}(dx^1)^2+e^{x^1}(dx^2)^2+e^{x^1+x^2}(dx^3)^2+(dx^4)^2.
\ee
The non-zero components (upto symmetry) of $R$, $S$, $\kappa$ and $P$ are given by
$$R_{2323}=-\frac{1}{2} e^{x^1+x^2}; \ \ \ S_{22}=\frac{1}{2}, \ \ \ S_{33}=\frac{e^{x^2}}{2}; \ \ \ \kappa = e^{-x^1};$$
$$P_{1221}=-\frac{e^{x^1}}{6}, \  \ -2P_{1331}= -P_{2323}= P_{2332}=\frac{1}{3} e^{x^1+x^2}, \ \ P_{2424}=\frac{1}{6}, \ \ P_{3434}=\frac{e^{x^2}}{6}.$$
From above we can easily check the following:\\
(i) $R\cdot R = 0$ and thus $R\cdot S = 0$, $P\cdot S =0$ and $R\cdot P =0$.\\
(ii) $R = S\wedge S$ and thus $Q(S,R) = 0$.\\
(iii) As here $R\cdot R = 0$ and $Q(S,R) = 0$ so $P\cdot R = 0$ and $P\cdot P = -\frac{1}{3} Q(S,P)$.\\
(iv) Although $P\cdot R =0$ and $P\cdot S =0$ but $P\cdot \mathcal R \ne 0$ and also $P\cdot \mathcal S \ne 0$.\\
%--------------------------------------------------------------------------
\textbf{Note:} This example ensures that on a semi-Riemannian manifold, the curvature conditions $P\cdot R = 0$ and $P\cdot \mathcal R = 0$ give different structures.\\
%=============================================================================================================
\textbf{Example 2:} Let $M_2$ be a $4$-dimensional connected semi-Riemannian manifold endowed with the semi-Riemannian metric
\be\label{met2}
ds^2 = (1+2e^{x^1})\left[(dx^1)^2+(dx^2)^2+(dx^3)^2+(dx^4)^2\right].
\ee
Then the non-zero components (upto symmetry) of $R$, $S$, $\kappa$ and $P$ are given by
$$R_{1212}= R_{1313}= R_{1414}= -\frac{e^{x^1}}{2 e^{x^1}+1}, \ \ R_{2323}= R_{2424}= R_{3434}= -\frac{e^{2 x^1}}{2 e^{x^1}+1};$$
$$S_{11}=\frac{3 e^{x^1}}{\left(2 e^{x^1}+1\right){}^2}, \ \ S_{22}= S_{33}= S_{44}=\frac{e^{x^1}}{2 e^{x^1}+1}; \ \ \kappa = \frac{6 e^{x^1} (1+e^{x^1})}{(1+2 e^{x^1})^3};$$
$$\frac{1}{2}P_{1221}= \frac{1}{2}P_{1331}= \frac{1}{2}P_{1441}= P_{2323}= -P_{2332}= P_{2424}= -P_{2442} = P_{3434}= -P_{3443}=-\frac{e^{2x^1}-e^{x^1}}{6 e^{x^1}+3}.$$
Using above we can easily calculate the non-zero components (upto symmetry) of $R\cdot R$, $Q(g,R)$, $Q(S,R)$ and $P\cdot R$ as follows:
$$R\cdot R_{122313}= R\cdot R_{122414}= -R\cdot R_{132312}= R\cdot R_{133414}= -R\cdot R_{142412}=-R\cdot R_{143413}=\frac{e^{2 x^1} \left(e^{x^1}-1\right)}{\left(2 e^{x^1}+1\right)^3};$$
$$Q(g,R)_{122313}= Q(g,R)_{122414}= -Q(g,R)_{132312}= Q(g,R)_{133414}= -Q(g,R)_{142412}= -Q(g,R)_{143413}$$$$=e^{x^1} \left(e^{x^1}-1\right);$$
$$Q(S,R)_{122313}= Q(S,R)_{122414}= -Q(S,R)_{132312}= Q(S,R)_{133414}= -Q(S,R)_{142412}= -Q(S,R)_{143413}$$$$=\frac{e^{2 x^1} \left(e^{x^1}-1\right)}{\left(2 e^{x^1}+1\right)^3};$$
$$P\cdot R_{122313}= P\cdot R_{122414}= -P\cdot R_{132312}= P\cdot R_{133414}= -P\cdot R_{142412}= -P\cdot R_{143413}=\frac{2 e^{2 x^1} \left(e^{x^1}-1\right)}{3 \left(2 e^{x^1}+1\right)^3}.$$
In view of above results we have the following pseudosymmetric type conditions on $M_2$:\\
(i) $R\cdot R = \frac{e^{x^1}}{\left(2 e^{x^1}+1\right)^3} Q(g,R) = Q(S,R)$ and hence $R\cdot S = P\cdot S = \frac{e^{x^1}}{\left(2 e^{x^1}+1\right)^3} Q(g,S)$\\
(ii) $P\cdot R = \frac{2 e^{x^1}}{3 \left(2 e^{x^1}+1\right)^3} Q(g,R) = \frac{2}{3}Q(S,R)$\\
(iii) $R\cdot R = L Q(g,R) + \left[1- L e^{-x^1} \left(2 e^{x^1}+1\right)^3\right]Q(S,R)$, $L$ being arbitrary scalar on $M$.\\
%
%====================================================Example 3================================================
%=============================================================================================================
\textbf{Example 3:} Let $M_3$ be a $5$-dimensional connected semi-Riemannian manifold endowed with the semi-Riemannian metric
\be\label{met3}
ds^2 = a(dx^1)^2+e^{2x^2}(x^4)^2(dx^2)^2+2e^{2x^2}dx^2 dx^3+e^{2x^2}(dx^4)^2+ (e^{2x^2}f) (dx^5)^2,
\ee
where $a$ is a positive constant and $f$ is a positive function of $x^2$.\\
%=============================================
The non-zero components (upto symmetry) of $R$, $S$, $\kappa$ and $P$ are given by
$$-R_{1212}=e^{2 x^1} (x^4)^2, \ \ -R_{1213}= -R_{1414}=e^{2 x^1}, \ \ -R_{1515}=f e^{2 x^1}, \ \ R_{2323}= -R_{2434}=\frac{e^{4 x^1}}{a},$$
$$-R_{2424}=\frac{e^{2 x^1} \left(a+e^{2 x^1} (x^4)^2\right)}{a}, \ \ -R_{2525}=\frac{e^{2 x^1} \left(2 a f f''-a \left(f'\right)^2+4 f^2 e^{2 x^1} (x^4)^2\right)}{4 a f},$$
$$-R_{2535}= -R_{4545}=\frac{f e^{4 x^1}}{a};$$
%---------------------------------------------
$$S_{11}=4, \ \ S_{22}=\frac{2 a f f''-a \left(f'\right)^2+4 a f^2+16 f^2 e^{2 x^1} (x^4)^2}{4 a f^2}, \ \ 
S_{23}= S_{44}=\frac{4 e^{2 x^1}}{a}, \ \ S_{55}=\frac{4 f e^{2 x^1}}{a};$$
%---------------------------------------------
$$\kappa = \frac{20}{a}; \ \ P_{1221}=-\frac{a \left(2 f f''-\left(f'\right)^2+4 f^2\right)}{16 f^2}, \ \ P_{2322}=\frac{e^{2 x^1} \left(2 f f''-\left(f'\right)^2+4 f^2\right)}{16 f^2},$$
$$P_{2424}=-\frac{e^{2 x^1} \left(-2 f f''+\left(f'\right)^2+12 f^2\right)}{16 f^2},$$
$$P_{2442}=e^{2 x^1}, \ \ P_{2525}=\frac{e^{2 x^1} \left(-6 f f''+3 \left(f'\right)^2+4 f^2\right)}{16 f}, \ \ P_{2552}=\frac{e^{2 x^1} \left(2 f f''-\left(f'\right)^2\right)}{4 f}.$$
%==============================================
The non-zero components (upto symmetry) of $R\cdot R$ and $Q(g,R)$ are given by
$$R\cdot R_{122414}=-R\cdot R_{142412}=e^{2 x^1}, \ \ R\cdot R_{122515}= -R\cdot R_{152512}=\frac{e^{2 x^1} \left(2 f f''-\left(f'\right)^2\right)}{4 f}$$
$$-2R\cdot R_{232424}= R\cdot R_{242423}=\frac{2 e^{4 x^1}}{a}, \ \ -R\cdot R_{232525}= \frac{1}{2}R\cdot R_{252523}= R\cdot R_{254524}=\frac{e^{4 x^1} \left(2 f f''-\left(f'\right)^2\right)}{4 a f}$$
$$R\cdot R_{242545}=\frac{e^{4 x^1} \left(-2 f f''+\left(f'\right)^2+4 f^2\right)}{4 a f}, \ \ R\cdot R_{244525}=-\frac{f e^{4 x^1}}{a},$$
%-----------------------------------------------
$$Q(g,R)_{122414}= -Q(g,R)_{142412}=a e^{2 x^1}, \ \ Q(g,R)_{244525}=-f e^{4 x^1},$$
$$Q(g,R)_{122515}= -Q(g,R)_{152512}=\frac{a e^{2 x^1} \left(2 f f''-\left(f'\right)^2\right)}{4 f},$$
$$-2Q(g,R)_{232424}= Q(g,R)_{242423}=2 e^{4 x^1}, \ \ Q(g,R)_{242545}=\frac{e^{4 x^1} \left(-2 f f''+4 f^2+\left(f'\right)^2\right)}{4 f},$$
$$-Q(g,R)_{232525}= \frac{1}{2}Q(g,R)_{252523}= Q(g,R)_{254524}=\frac{e^{4 x^1} \left(2 f f''-\left(f'\right)^2\right)}{4 f}.$$
%==================================================
It can be now easily checked that $P\cdot \mathcal S = 0$. From the components of above tensors, we see that the semi-Riemannian manifold $M_3$ satisfies the following geometric structures:\\
(i) It is a Venzi's $P$-space satisfying
$$\Pi(X_1)P(X_2,X_3,X,Y)+\Pi(X_2)P(X_3,X_1,X,Y)+\Pi(X_3)P(X_1,X_2,X,Y)=0$$
for the null 1-form $\Pi = \{0,1,0,0,0\}$. Here $W\cdot W = 0$ and hence improperly satisfies the condition $W\cdot W = L Q(n S-\kappa g, W)$, which supports the Theorem \ref{thm4.3}.\\
(ii) $R\cdot R = \frac{1}{a}Q(g,R)$, i.e., pseudosymmetric manifold of constant type and thus Ricci, conformally, projectively, concircularly and conharmonically pseudosymmetric manifoold of constant type.\\
(iii) $P\cdot S = \frac{1}{a}Q(g,S) \ne 0$ but $P\cdot\mathcal S=0$.\\
%
%====================================================Example 4================================================
%=============================================================================================================
\textbf{Example 4:} Let $M_4$ be a $4$-dimensional connected semi-Riemannian manifold endowed with the semi-Riemannian metric
\be\label{met4}
ds^2 = x^1 x^3 (dx^1)^2+2 dx^1 dx^2 + (2+dx^1)^2 dx^3+ (x^1)^3(dx^4)^2.
\ee
%=============================================
The non-zero components (upto symmetry) of $R$, $S$, $\kappa$ and $P$ are given by
$$R_{1414}=-\frac{3 x^1}{4}; \ \ \ S_{11}=\frac{3}{4 (x^1)^2}; \ \ \kappa =0;$$
$$P_{1211}=\frac{1}{4 (x^1)^2}, \ \ P_{1313}=\frac{\left(x^1+2\right)^2}{4 (x^1)^2}, \ \ P_{1441}= -\frac{3}{2}P_{1414}=\frac{3 x^1}{4}.$$
From the components of above tensors it is easy to check that $R\cdot R =0$, $R\cdot S = 0$, $Q(S,R)=0$, $P\cdot R = 0$ and also $P\cdot \mathcal S = 0$. According to the Theorem \ref{thm4.5}, if the manifold is of non-constant scalar curvature, then $P\cdot S = 0$ and $P\cdot \mathcal S =0$ holds simultaneously, then the manifold is Einstein. This example also supports the Theorem \ref{thm4.8}.
%%%%%%%%%%%%%%%%%%%%%%%%%%%%%%%%%%%%%%%%%%%%%%%%%%%%%%%%%%%%%%%%%%%%%%%%%%%%%%%%%%%%%%%%%%%%%%%%%%%%%%%%%%%%%%
%                                                  Conclusions
%%%%%%%%%%%%%%%%%%%%%%%%%%%%%%%%%%%%%%%%%%%%%%%%%%%%%%%%%%%%%%%%%%%%%%%%%%%%%%%%%%%%%%%%%%%%%%%%%%%%%%%%%%%%%%%%
\section{\bf{Conclusions}}
%=========================
In the present paper we study the basic properties of the projective operator and calculate the necessary and sufficient conditions for a semi-Riemannian manifold to satisfy some Walker type conditions. It is shown that the projective operator commutes with contraction if and only if the manifold is Einstein. A necessary and sufficient condition for a semisymmetric (resp., pseudosymmetric) manifold due to projective curvature tensor to be a Ricci semisymmetric (resp., pseudosymmetric) manifold due to projective curvature tensor is presented. It is also shown that $P \cdot \mathcal R = 0$ and $P\cdot R=0$ (resp., $P \cdot \mathcal R = L Q(g,\mathcal R)$ and $P\cdot R= L Q(g,R)$) give different structures. We have evaluated some pseudosymmetric type condition due to projective curvature tensor under certain condition and showed that a $P$-space by Venzi is either $W$-space by Venzi or a manifold of constant curvature. We obtain the curvature properties of of various semisymmetric type and pseudosymmetric type curvature restricted geometric structures due to projective curvature tensor, such as 
(i) $P\cdot R = 0,$ 
(ii) $P\cdot R = L Q(g,R)$, 
(iii) $P\cdot S = 0,$ 
(iv) $P\cdot S = L Q(g, S)$, 
(v) $P\cdot \mathcal S = 0,$ 
(vi) $P\cdot \mathcal S = L Q(g,\mathcal S)$, 
(vii) $P\cdot \mathcal R = 0$, 
(viii) $P\cdot \mathcal R = L Q(g,\mathcal R)$, 
(ix) $P\cdot \mathcal R = L Q(S,\mathcal R)$, 
(x) $P\cdot P = 0,$ 
(xi) $P\cdot P = L Q(g,P)$,  
(xii) $P\cdot P = L Q(S,P)$, 
(xiii) $P\cdot \mathcal P = 0$, 
(xiv) $P\cdot \mathcal P = L Q(g,\mathcal P)$ and 
(xv) $P\cdot \mathcal P = L Q(S,\mathcal P)$ on a Riemannian as well as semi-Riemannian manifold. It is shown that a Riemannian manifold $M$ with one of the curvature condition (v)-(xv) reduce to a Einstein manifold and hence manifold of constant curvature if $M$ is generalized Roter type.
%===============================================================================================================
%\textbf{Acknowledgment:} The second named author gratefully acknowledges to CSIR, New Delhi (File No. 09/025 (0194)/2010-EMR-I) for the financial assistance. All the algebraic computations of Section \ref{exam} are performed by a program in Wolfram Mathematica.
%%%%%%%%%%%%%%%%%%%%%%%%%%%%%%%%%%%%%%%%%%%%%%%%%%%%%%%%%%%%%%%%%%%%%%%%%%%%%%%%%%%%%%%%%%%%%%%%%%%%%%%%%%%%%%%%%
%%%%%%%%%%%%%%%%%%%%%%%%%%%%%%%%%%%%%%%%%%%%%%%%%%%%%%%%%%%%%%%%%%%%%%%%%%%%%%%%%%%%%%%%%%%%%%%%%%%%%%%%%%%%%%%%%

%%%%%%%%%%%%%%%%%%%%%

%%%%%%%%%%%%%%%%%%%%%%%%%%%%%%%%%%%%%%%%%%%%%%%%%%%%%%%%%%%%%%%%%%%%%%%%%%%%%%%%%%%%%%%%%%%%%%%%%%%%
\end{document}